\documentclass[%
 aip,
 amsmath,amssymb,
 reprint,%
]{revtex4-1}

\usepackage{graphicx}% Include figure files
\usepackage{dcolumn}% Align table columns on decimal point
\usepackage{bm}% bold math
\usepackage{mathrsfs}
\usepackage[utf8]{inputenc}
\usepackage[T1]{fontenc}
\usepackage{mathptmx}
\usepackage{etoolbox}
\usepackage[title]{appendix}
\newtheorem{thm}{Theorem}[section]
\newtheorem{remark}{Remark}[section]
\newtheorem{pro}{Proof}
\newtheorem{proposition}{Proposition}[section]
\newtheorem{lemma}{Lemma}[section]
\newtheorem{definition}{Definition}[section]

\numberwithin{equation}{section}
\makeatletter
\def\@email#1#2{%
 \endgroup
 \patchcmd{\titleblock@produce}
  {\frontmatter@RRAPformat}
  {\frontmatter@RRAPformat{\produce@RRAP{*#1\href{mailto:#2}{#2}}}\frontmatter@RRAPformat}
  {}{}
}%
\makeatother
\begin{document}

\preprint{AIP/123-QED}

\title{The Cauchy problem of non-local space-time reaction-diffusion\\equation
involving fractional $p$-Laplacian }
% Force line breaks with \\
\author{Fei. Gao}
\email{gaof@whut.edu.cn}
 \altaffiliation{Author to whom correspondence should be addressed: gaof@whut.edu.cn}

\author{Hui Zhan}
\altaffiliation{Electronic mail:2432593867@qq.com}
\affiliation{%
Department of Mathematics and Center for Mathematical Sciences,\\ Wuhan University of Technology, Wuhan, 430070, China%\\This line break forced% with \\
}%

\date{\today}% It is always \today, today,

\begin{abstract}
For the  non-local space-time reaction-diffusion equation
involving fractional $p$-Laplacian \begin{equation*}
    \begin{cases}
          \frac{\partial^{\alpha }u}{\partial t^{\alpha }}+(-\Delta)_{p}^{s} u=\mu u^{2}(1-kJ*u)-\gamma u,&(x,t)\in\mathbb{R}^{N}\times(0,T)\\
u(x,0)=u_{0}(x),& x\in\mathbb{R}^{N}
    \end{cases}
\end{equation*} $\mu>0 ,k>0,\gamma\geq 1,\alpha\in(0,1),s\in(0,1),1<p$, we consider for $N\leq2$ the problem of finding
a global boundedness of the weak solution by virtue of Gagliardo-Nirenberg
inequality and fractional Duhamel's formula. Moreover, we prove such weak
solution converge to $0$ exponentially or locally uniformly as $t \rightarrow
\infty$ for small $\mu$ values  with the comparison principle and local
Lyapunov type functional. In those cases the problem is reduced to   fractional
$p$-Laplacian equation in the non-local reaction-diffusion range  which is
treated with the symmetry and other properties of the kernel of
$(-\Delta)_{p}^{s}$. Finally, a key element in our construction is a proof of
global bounded weak solution with the fractional nonlinear diffusion terms
$(-\Delta)_{p}^{s}u^{m}(2-\frac{2}{N}<m\leq 3,1<p<\frac{4}{3})$  by using Moser iteration and
fractional differential inequality. 
\end{abstract}

\maketitle
\section{Introduction}
In recent years, the study of differential equations using non-local fractional operators has attracted a lot of interest \cite{borikhanov2022qualitative}. The space-time fractional equations could be applied to a wide range of applications such as: continuum mechanics, phase transition phenomena, population dynamics, image process, game theory and Lévy processes, see \cite{2021Backward,Sal2014Simultaneous,2013Aqqq,2018An} and the references therein. In particular,  by \cite{2018Regularity}, we know that the fractional Laplacian operator is defined for  functions $u(x),x\in \mathbb{R}^{N}$, as 
\begin{align*}
    (-\Delta)^{s} u(x)=C_{n,s}P.V.\int _{\mathbb{R}^{N}}\frac{u(x)-u(y)}{\left | x-y \right |^{n+2s}}dy,
\end{align*}	
where $s\in(0,1),C_{n,s}=\frac{4^{s }s \Gamma (\frac{n}{2}+s)}{\Gamma (1-s )\pi ^{\frac{n}{2}}},$ and $P.V.$ is the principal value of Cauchy. In the special case when $p=2$, Then problem \eqref{1} reduces to the space-time fractional
non-local reaction-diffusion equation \cite{hui2022global}.

 We are concerned here with the nonlinear version given by the fractional $p$-Lapalcian operator $(-\Delta)_{p}^{s}$ defined by Definition \ref{opliyredfv}. To be precise, it can be called the $s$-fractional $p$-Laplacian operator \cite{borikhanov2022qualitative}. Following our previous papers \cite{2021zhan,hui2022global}, we will continue the study of the cauchy problem of weak solution,i.e., the space-time reaction-diffusion
equation involving fractional p-Laplacian
\begin{align} &\frac{\partial^{\alpha }u}{\partial t^{\alpha
}}+(-\Delta)_{p}^{s} u=\mu u^{2}(1-kJ*u)-\gamma
u, \label{1}\\
&u(x,0)=u_{0}(x),\quad x\in\mathbb{R}^{N} , \label{2}
	\end{align}
where $\mu ,k>0,\gamma\geq 1,0< \alpha,s<1, N\leq 2,1<p$.  And furthermore,
for $\alpha=1,s=1$ we recover the classical non-local reaction-diffusion
equation, whose theory is well known,cf \cite{0Global}. We also assume that we
are given the  competition kernel $J(x)$ with \begin{equation}\label{3}
  0\leq J \in L^{1}(\mathbb{R}^{N}),\, \int_{\mathbb{R}^{N}}J(x)dx=1,\, \underset{B(0,\delta_{0})}{inf}J>\eta,
  \end{equation}
for some $\delta_{0}>0,\eta>0$ and where $B(0,\delta_{0})=(-\delta_{0},\delta_{0})^{N}$ and
 $$J*u(x,t)=\int_{\mathbb{R}^{N}}J(x-y)u(y,t)dy.$$
 The related equations of this model are currently of great interest to researchers. There, the case $p\geq 2$, usually called nonlinear time-space fractional diffusion equation, has been treated in \cite{borikhanov2022qualitative}. We pursue in this paper the analysis of such a space-time diffusion  equation involving fractional $p$-Laplacian devote much attention to settle Cauchy problem of equation \eqref{1}-\eqref{2}, like blow-up, global boundedness and asymptotic behavior. Moreover, denoting
$J=k=\mu=\gamma=1$, we  verify  existence and global boundedness of weak solution for problem
\eqref{114}-\eqref{1.1.5}  in  fractional Sobolev space.

\textbf{Outline of result.}
As a starting novelty, the present paper recall some necessary definitions and useful properties of the fractional Sobolev space. See whole details in Section \ref{wedvfg}. Note that for $1<p$, we prove the blow-up, global boundedness and asymptotic behavior in equation \eqref{1}-\eqref{2}. However, in the range $1<p<\frac{4}{3}$, we need the extra condition, $sp<1$, that seems to play a major role in the proof for Theorem \ref{theorem6}. Indeed, in order to proof the existence of equation \eqref{114}-\eqref{1.1.5}, we introduces a basic tool called the weighted $L^{1}$ estimate that will play an important role in the existence theory for general classes of data in Appendix \ref{zhanhui22} .

In section \ref{90oplkjh}, we aim to proof the global boundedness for equation
\eqref{1}-\eqref{2}. The proof consists of several steps. First,  by reference \cite{borikhanov2022qualitative}, we obtain the global existence of weak solution.
Further, we prove the blow-up of the problem by reference
\cite{2020Global11,2014Blowing}. Next, two key steps in our proofs are briefly
listed in the following. Let $\delta_{0},\mu>0$ be as in \eqref{3} and $x\in
\mathbb{R}^{N}$ fixed. For any $0<\delta\leq \frac{1}{2}\delta_{0}$, using the
symmetry of the kernel of $u\in W^{s,p}(B(x,\delta))$ and the observation:
$\int_{B(x,\delta )} u\left ( -\Delta  \right )_{p}^{s}udy=\frac{1}{2}\left [ u
\right ]_{W^{s,p}(B(x,\delta ))}^{p}$, we obtain a key differential inequality
on $\int _{B(x,\delta )}u^{2}dy$ (see \eqref{11}) \begin{align*}
&\frac{\partial^{\alpha } }{\partial t^{\alpha }}\int_{B(x,\delta
)}u^{2}dy+\left [ u \right ]_{W^{s,p}(B(x,\delta ))}^{p}+2\gamma\int_{B(x,\delta )}u^{2}dy\\ &\leq
2\mu\int_{B(x,\delta )}u^{3}dy-2\mu\eta k\int_{B(x,\delta
)}u^{3}dy\int_{B(x,\delta )}udy. 
\end{align*}
with further controlling of the term $\mu\int_{B(x,\delta )}u^{3}dy$ by the
Gagliardo-Nirenberg inequality. Then we can obtain the uniform boundedness of
$\int_{B(x,\delta )}u^{2}dy$ by using Sobolev embedding inequality and
fractional differential equation. Based on the these, by applying fractional
Duhamel's formula representation in time spitted intervals $[0,T)$, for
example, let $f(u)=\mu u^{2}(1-kJ*u)-\gamma u$, we have the solution
 $$u(t)=\mathcal{S}_{\alpha }(t)u_{0}+\int_{0}^{t}(t-s)^{\alpha -1}\mathcal{K} _{\alpha }(t-s)f(u(s))ds,$$
where $\mathcal{S}_{\alpha },\mathcal{K} _{\alpha }$ is defined by Lemma
\ref{lem3}. Finally, we can obtain the global uniform boundedness of $u$ by a
series of careful analyses of the heat kernel.

In section \ref{0kjhgfrt}, in order to study the asymptotic behavior of weak
solutions, we  use the global boundedness of $u$ of Section  \ref{90oplkjh} to
obtain $0\leq u(x,t)<a$ for any $(x,t)\in \mathbb{R}^{N}\times [0,+\infty )$.
Then, we can derive a Lyapunov type functional on $B(x,\delta)$ in Proposition
\ref{proposition3.1}
 \begin{equation}\label{oplkjyt}
     \frac{\partial^{\alpha} H(x,t)}{\partial t^{\alpha}}\leq -(-\Delta)_{p}^{s} H(x,t)+\int _{B(x,\delta)}\left | \nabla u(y,t) \right |^{2}dy-D(x,t)
   \end{equation}
with
   $$D(x,t)=\frac{1}{2}(A-a)\mu k\int _{B(x,\delta)}u^{2}(y,t)dy.$$
This Lyapunov function plays key roles in studying the long time behavior.
Thus, with the help of comparison principle of equation \eqref{1}-\eqref{2} and
using fractional Duhamel's formula to \eqref{oplkjyt}, we can proceed to deduce
the long time behavior of the solution.

Section \ref{9o0oiyhnn} studies the global boundedness of equation
\eqref{114}-\eqref{1.1.5} relying on fractional differential inequality,
Gagliardo-Nirenberg inequality and other common inequalities. Before the proof
we used weighted $L^{1}$ estimates  to prove the existence and uniqueness of
the weak solution in Appendix \ref{zhanhui22}-\ref{28uijybb}. Firstly, due to
Sobolev inequality, one has a key formula that \begin{align}\label{olpkhftg}
\nonumber    &\frac{1}{k}(_{0}^{C}\textrm{D}_{t}^{\alpha }\int
_{\mathbb{R}^{N}}u^{k}dx)+\frac{c_{3}S}{2}\left ( \int _{\mathbb{R}^{N}} \left
| u^{\frac{m(p-1)+k-1}{p}}\right |^{p_{s }^{*}}dx\right
)^{\frac{p}{p_{s}^{*}}}\\
    &\leq\int _{\mathbb{R}^{N}}u^{k+1}dx(1-\int _{\mathbb{R}^{N}}udx)
-\int _{\mathbb{R}^{N}}u^{k}dx. \end{align} And by using fractional
differential inequality (see Lemma \ref{lemma444}-\ref{4554}), we obtains the
$L^{r}$ estimates for equation \eqref{114}-\eqref{1.1.5}. Next, in order to
extend $L^{\infty}$ estimation to $L^{\infty}$ estimation, by taking
$k=q_{k}:=2^{k}+2$ in \eqref{olpkhftg}, we obtain \eqref{4.2.1}. Furthmore,
through Moser iteration and $L^{r}$ estimation, the global boundedness of the
solution for \eqref{114}-\eqref{1.1.5} is obtained in fractional Sobolev space.

As already mentioned above,  we collect some useful lemma about  fractional derivative , fractional Duhamel's formula and blow-up for equation \eqref{1}-\eqref{2} in Appendix \ref{0opliutv}-\ref{0po9754}. Furthmore, The Appendix
\ref{i9ojjg} gives some useful results for proof of theorem. Finally,  based on
weighted $L^{1}$ estimates in Appendix \ref{zhanhui22},  we study the Appendix \ref{28uijybb} to obtain
the existence of the weak solution for \eqref{114}-\eqref{1.1.5}.

\textbf{NOTATIONS:}  In the whole paper we fix Caputo derivative $\partial
_{t}^{\alpha },\alpha\in (0,1)$ models memory effects in time. And we also
replace fractional $p$-Laplacian $-(-\Delta )_{p}^{s},s\in(0,1),1<p$ with  Neumann Laplacian operator $\Delta$ in
section \ref{90oplkjh}-\ref{0kjhgfrt}. In contrast, under the conditions of
$1<p<\frac{4}{3}$, the nonlinear fractional  term $-(-\Delta)_{p}^{s} u^{m}$
satisfies $2-\frac{2}{N}<m\leq 3$.

\section{ Preliminaries}\label{wedvfg}
\textbf{The fractional Sobolev space.}\quad According to \cite{mazon2016fractional}, let
$\Omega$ be an open set in $\mathbb{R}^{N}$. For any $max\left \{
\frac{2N}{N+2s},1 \right \}<p$ and any $0<s<1$, let us denote by
	$$\left [ u \right ]_{W^{s,p}(\Omega )}=\left ( \int _{\Omega }\int _{\Omega }\frac{\left | u(x)-u(y) \right |^{p}}{\left | x-y \right |^{N+sp}}dxdy \right )^{\frac{1}{p}},$$
the $(s,p)$-Gagliardo seminorm of a measurable function $u$ in $\Omega$. We
consider the fractional Sobolev space $$W^{s,p}(\Omega )=\left \{ u\in
L^{p}(\Omega ):\left [ u \right ]_{W^{s,p}(\Omega )}<\infty  \right \},$$ which
is a Banach space with respect to the norm $$\left \| u \right
\|_{W^{s,p}(\Omega )}:=\left [ u \right ]_{W^{s,p}(\Omega )}+\left \| u \right
\|_{L^{p}(\Omega )}.$$

\begin{definition}\textsuperscript{\cite{Kilbas2006}}

\textbf{(i)}\quad Assume that $X$ is a Banach space and let $u:\left [ 0,T
\right ]\rightarrow X$. The Caputo fractional derivative operators of $u$ is
defined by \begin{eqnarray} _{0}^{C}\textrm{D}_{t}^{\alpha
}u(t)&=&\frac{1}{\Gamma (1-\alpha )}\int_{0}^{t}(t-s)^{-\alpha
}\frac{d}{ds}u(s)ds,\label{32}\\ \nonumber_{t}^{C}\textrm{D}_{T}^{\alpha
}u(t)&=&\frac{-1}{\Gamma (1-\alpha )}\int_{t}^{T}(s-T)^{-\alpha
}\frac{d}{ds}u(s)ds, \end{eqnarray}
 where $\Gamma (1-\alpha )$ is the Gamma function. The above integrals are called the left-sided and the right-sided the Caputo fractional derivatives.

 \textbf{(ii)}\quad For $u:[ 0,\infty )\rightarrow \mathbb{R}^{N}$, the left Caputo fractional derivative with respect to time $t$ of $u$ is defined by
 \begin{equation}\label{44}
  \partial _{t}^{\alpha }u(x,t)=\frac{1}{\Gamma (1-\alpha )}\int_{0}^{t}\frac{\partial u}{\partial s}(x,s)(t-s)^{-\alpha }ds,\quad 0<\alpha<1\quad
\end{equation} \end{definition}

\begin{definition}\textsuperscript{\cite{Kilbas2006}}
 The left and right Rieman-Liouville fractional integrals of order $\alpha\in (0,1)$ for an integrable function $u(t)$ are given by
 \begin{equation*}
     _{0}\textrm{I}_{t}^{\alpha }u(t)=\frac{1}{\Gamma (\alpha )}\int_{0}^{t}(t-s)^{\alpha -1}u(s)ds,\quad t\in (0,T],
 \end{equation*}
 and
 \begin{equation*}
      _{t}\textrm{I}_{T}^{\alpha }u(t)=\frac{1}{\Gamma (\alpha )}\int_{t}^{T}(s-t)^{\alpha -1}u(s)ds,\quad t\in [0,T).
 \end{equation*}
\end{definition}

\begin{definition}\textsuperscript{\cite{Kilbas2006}}\label{proposition 2.3}
  The Mittag-Leffler function $E_{\alpha}(z),E_{\alpha ,\beta }(z)$ is defined as
   $$E_{\alpha}(z)=\sum_{k=0}^{\infty }\frac{z^{k}}{\Gamma (\alpha k+1)},\quad E_{\alpha ,\beta }(z)=\sum_{k=0}^{\infty }\frac{z^{k}}{\Gamma (\alpha k+\beta )},$$
  where$\beta, z\in \mathbb{C},R(\alpha)>0,\mathbb{C}$ denote the complex plane.
  \end{definition}
\begin{definition}\textsuperscript{\cite{mazon2016fractional}}\label{opliyredfv}
   The fractional $p$-Laplacian operator for $s\in(0,1),p>1$ and $u\in W^{s,p}(\mathbb{R}^{N})$, is defined by
\begin{align}\label{4} 
\nonumber   (-\Delta)_{p}^{s} u&=P.V.\int
_{\mathbb{R}^{N}}\frac{\left | u(x)-u(y) \right |^{p-2}(u(x)-u(y))}{\left | x-y
\right |^{N+sp}}dy\\ &=\lim_{\varepsilon \downarrow 0}\int
_{\mathbb{R}^{N}\setminus B_{\varepsilon }(x)}\frac{\left | u(x)-u(y) \right
|^{p-2}(u(x)-u(y))}{\left | x-y \right |^{N+sp}}dy .
\end{align}
where $P.V.$ is the principal value of Cauchy.
\end{definition}
\begin{thm}\textsuperscript{\cite{borikhanov2022qualitative}}\label{theorem1}
 (Existence)\,Let $u_{0}\in W_{0}^{s,p}(\Omega ),u_{0}\geq 0,sp<N$, then there exists $T>0$ such that the problem \eqref{1}-\eqref{2} has a local real-valued unique weak solution $u\in \Pi $ on $(0,T)$, where $\Pi $ is defined by
    $$\Pi =\left \{ u,\partial _{t}^{\alpha }u\in W^{s,p}(\Omega ;L^{\infty }(0,T))\cap L^{2}(\Omega;L^{\infty }(0,T) ) \right \}.$$
\end{thm}

\section{Global boundedness of weak solutions }\label{90oplkjh}
\begin{thm}\label{theorem4}
   Assume \eqref{3} holds and  $0<u_{0}\in X_{0},T>0$. When $k>k^{*}$, denoting  $k^{*}=0$ for $N=1$ and $k^{*}=(\mu C_{GN}^{2}+1)\eta ^{-1}$ for $N=2$, where $C_{GN}$  is the constant appears in Gagliardo-Nirenberg inequality in Lemma \ref{lemma2.1}, there exists a non-negative weak solution to equation \eqref{1}-\eqref{2} and it is globally bounded in time, that is, there exists
  \begin{equation}
 \nonumber   K=\begin{cases}
 K(\Arrowvert u_{0}\Arrowvert_{L^{\infty}(\mathbb{R}^{N})},\mu,\eta,k,C_{GN}),&  N=1,\\
 K(\Arrowvert u_{0}\Arrowvert_{L^{\infty}(\mathbb{R}^{N})},\mu,C_{GN}),&  N=2,
      \end{cases}
  \end{equation}
  such that
\begin{equation*}
     0\leq u(x,t)\leq K,\qquad\forall (x,t)\in \mathbb{R} ^{N}\times (0,T).
  \end{equation*}
\end{thm}

\begin{pro}

 First, we denote $$B(x_{0},\delta  ):=\left \{ x \triangleq
(x_{1},\cdots ,x_{N})\in \mathbb{R}^{N}|\left | x_{i}-x_{i}^{0}  \right |\leq
\delta,1\leq i\leq N\right \}$$ where $x_{0}\triangleq(x_{1}^{0},\cdots
,x_{N}^{0})\in \mathbb{R}^{N}$ and $0<\delta\leq \frac{1}{2}\delta_{0}$. And we
choose $\varphi_{\varepsilon}(\cdot)\in C_{0}^{\infty}(B(x,\delta ))$, and
$\varphi_{\varepsilon}(\cdot)\to 1$ locally uniformly in $B(x,\delta )$ as
$\varepsilon\to 0$. For any $x \in \mathbb{R}^{N}$, we multiply \eqref{1} by
$2u\varphi_{\varepsilon}$ and integrate by parts over $B(x,\delta )$, we obtain

\begin{align*}
&\int_{B(x,\delta
)}2u\varphi_{\varepsilon}\frac{\partial^{\alpha } u}{\partial t^{\alpha
}}dy+\int_{B(x,\delta )}2u\varphi_{\varepsilon}(-\Delta)_{p}^{s}udy\\
&=\int_{B(x,\delta )}2u^{3}\mu\varphi_{\varepsilon}(1-kJ*u)dy-2\gamma\int _{B(x,\delta )}u^{2}\varphi_{\varepsilon} dy. 
\end{align*}
By Lemma \ref{lemma23}, then we can get $$\int_{B(x,\delta
)}2u\varphi_{\varepsilon}\frac{\partial^{\alpha } u}{\partial t^{\alpha
}}dy\geq\frac{\partial^{\alpha }}{\partial t^{\alpha }}\int_{B(x,\delta
)}u^{2}\varphi_{\varepsilon}dy.$$ And we assume that the integral in the
definition of $(-\Delta)_{p}^{s}$ exists, then for the symmetry of the kernel
of $u\in W^{s,p}(B(x,\delta ))$, we obtain the following formula for the
integration by parts 
\begin{align}\label{7} \nonumber  & \int_{B(x,\delta )}
u\left ( -\Delta  \right )_{p}^{s}udy\\
\nonumber&=\int _{B(x,\delta )}\int _{B(x,\delta
)}\frac{\left | u(x,t)-u(y,t) \right |^{p-2}(u(x,t)-u(y,t))}{\left | x-y \right
|^{N+sp}}\\
\nonumber&\qquad u(x,t)dxdy\\
\nonumber&=\int _{B(x,\delta )}\int _{B(x,\delta
)}\frac{\left | u(x,t)-u(y,t) \right |^{p-2}(u(y,t)-u(x,t))}{\left | x-y \right
|^{N+sp}}\\
\nonumber&\qquad u(y,t)dxdy\\ 
\nonumber&=-\int _{B(x,\delta )}\int _{B(x,\delta
)}\frac{\left | u(x,t)-u(y,t) \right |^{p-2}(u(x,t)-u(y,t))}{\left | x-y \right|^{N+sp}}\\
\nonumber&\qquad  u(y,t)dxdy\\ 
\nonumber&=\frac{1}{2}\int _{B(x,\delta )}\int
_{B(x,\delta )}\frac{\left | u(x,t)-u(y,t) \right |^{p}}{\left | x-y \right
|^{N+sp}}dxdy\\ 
&=\frac{1}{2}\left [ u \right ]_{W^{s,p}(B(x,\delta ))}^{p}.
\end{align} 
By \eqref{7} and $\varepsilon\to 0$, we have 
\begin{align*}
&\frac{\partial^{\alpha } }{\partial t^{\alpha }}\int_{B(x,\delta
)}u^{2}dy+\left [ u \right ]_{W^{s,p}(B(x,\delta ))}^{p}\\
&\leq
2\mu\int_{B(x,\delta )}u^{3}(1-kJ*u)dy-2\gamma\int_{B(x,\delta )}u^{2}dy,
\end{align*} 
where we applied the following facts that for any $y\in
B(x,\delta)$, then $J(z-y)\geq \eta$ and for any $z\in
B(y,2\delta),B(x,\delta)\subset B(y,2\delta)$, then \begin{equation*}
  J*u(y,t)\geq \eta \int _{B(y,2\delta )}u(z,t)dz\geq\eta\int_{B(x,\delta )}u(y,t)dy.
\end{equation*} 
Therefore, we get 
\begin{align}\label{11}
\nonumber&\frac{\partial^{\alpha } }{\partial t^{\alpha }}\int_{B(x,\delta
)}u^{2}dy+\left [ u \right ]_{W^{s,p}(B(x,\delta ))}^{p}\\ \nonumber&\leq
2\mu\int_{B(x,\delta )}u^{3}dy-2\mu\eta k\int_{B(x,\delta
)}u^{3}dy\int_{B(x,\delta )}udy\\
&\qquad-2\gamma\int_{B(x,\delta )}u^{2}dy.
\end{align}
Next, in order to estimate the term $\int_{\mathbb{R}^{N}}u^{3}dy$, by using
Gagliardo-Nirenberg inequality in Lemma \ref{lemma2.1}, we  combine
\eqref{22}-\eqref{25} to obtain
\begin{equation}
  \begin{aligned}
  &2\mu \int _{B(x,\delta )}u^{3}dy\\
  &\leq 2\int _{B(x,\delta )}\left | \nabla u \right |^{2}dy+2\mu(1+\mu^{\frac{N}{4-N}}C_{GN}^{\frac{4}{4-N}}(N,\delta))\\
  &\quad(\int _{B(x,\delta )}u^{3}dy\int _{B(x,\delta )}udy)^{\frac{6-N}{2(4-N)}}+2\mu C_{GN}^{\frac{2(6-N)}{N}}(N,\delta).\label{26}
  \end{aligned}
\end{equation} From this, we will consider the following $N=1$ and $N=2$ cases
respectively.

\textbf{Case 1.}  $N=1$. From \eqref{26}, by Young's inequality, let
$\varepsilon =\eta k,p=\frac{6}{5},q=6$, we get 
  \begin{align}\label{27}
 \nonumber  &2\mu \int _{B(x,\delta )}u^{3}dy\\
\nonumber   &\leq 2\int _{B(x,\delta )}\left | \nabla u \right |^{2}dy
+2\mu \eta k\int _{B(x,\delta )}u^{3}dy\int _{B(x,\delta )}udy\\ 
&\quad+2\mu
(\mu^{\frac{1}{3}}C_{GN}^{\frac{4}{3}}(1,\delta)+1)^{6}(\eta k)^{-5}+2\eta
C_{GN}^{10}(1,\delta).
  \end{align}
Bringing \eqref{27} into \eqref{11}, we have
   \begin{align}\label{28}
\nonumber  &\frac{\partial^{\alpha } }{\partial t^{\alpha }}\int_{B(x,\delta )}u^{2}dy+\left [ u \right ]_{W^{s,p}(B(x,\delta ))}^{p}+2\gamma\int _{B(x,\delta )}u^{2}dy\\
\nonumber&\leq 2\int _{B(x,\delta )}\left | \nabla u \right |^{2}dy\\
&\quad+2\mu
[(\mu^{\frac{1}{3}}C_{GN}^{\frac{4}{3}}(1,\delta))^{6}(\eta
k)^{-5}+C_{GN}^{10}(1,\delta)].
  \end{align}
 According to Lemma\ref{lemma99}, we know that there exists an embedding constant $C_{2}>0$ such that
 $$\left \| \nabla u \right \|_{L^{p}(B(x,\delta ))}^{2}\geq C_{1}\left \|  u\right \|_{L^{s_{2}}(B(x,\delta ))}^{2}$$
where $s_{2}\geq p$ will be determined later. Here we set $s_{2}=p=2$,then
there is
  $$\left \| \nabla u \right \|_{L^{2}(B(x,\delta ))}^{2}\geq C_{2}\left \|  u\right \|_{L^{2}(B(x,\delta ))}^{2}.$$
By \eqref{55yy}, we have \begin{equation*}
    \left [ u \right ]_{W^{s,p}(B(x,\delta ))}^{p}\geq C_{*}^{-p}\left \| u \right \|_{L^{2}(B(x,\delta ))}^{p}.
\end{equation*} Then, we let $$ Q_{1}=2\mu \left ( (\mu
^{\frac{1}{3}}C_{GN}^{\frac{4}{3}}(1,\delta)+1)^{6}(\eta
k)^{-5}+C_{GN}^{10}(1,\delta) \right ).$$ And the equation \eqref{28} is
configured above, we can get
\begin{align*}
     &\frac{\partial^{\alpha} }{\partial t^{\alpha}}\int _{B(x,\delta )}u^{2}dy+C_{*}^{-p}\left \| u \right \|_{L^{2}(B(x,\delta ))}^{p}+2\gamma\int _{B(x,\delta )}u^{2}dy\\
     &\leq C_{2}\int _{B(x,\delta )}u^{2}dy+Q_{1}.
\end{align*} Denote $w(t)=\int_{B(x,\delta )}u^{2}dy$, the solution of the
following  fractional differential equation
 \begin{equation*}
   \begin{cases}
    _{0}^{C}\textrm{D}_{t}^{\alpha }w(t)+(2\gamma-C_{2})w(t)+C_{*}^{-p}w^{\frac{p}{2}}(t)=Q_{1},\\
w(0)= 2\delta\left \| u_{0} \right \|_{L^{\infty }(\mathbb{R}^{N})}^{2},
   \end{cases}
 \end{equation*}
for any $ (x,t)\in \mathbb{R}\times [0,T_{max}]$.  From Lemma \ref{7700} and
Lemma \ref{lemma288}, we have
 $$
 \begin{aligned}
 &\int _{B(x,\delta )}u^{2}dy=\left \| u \right \|_{L^{2}(B(x,\delta ))}^{2}=  w(t)\\
&\leqslant w(0)+\left [
\lambda^{*}w(0)+\varepsilon^{\frac{2}{2-p}}C_{*}^{\frac{2p}{p-2}} +Q_{1}\right
]\\
&\quad\int_{0}^{t}(t-s)^{\alpha-1}E_{\alpha ,\alpha }(\lambda^{*}(t-s)^{\alpha})ds\\
&\leqslant2\delta \left \| u_{0} \right \|_{L^{\infty
}(\mathbb{R}^{N})}^{2}+\left [ 2\delta\lambda^{*}\left \| u_{0} \right\|_{L^{\infty
}(\mathbb{R}^{N})}^{2}+\varepsilon^{\frac{2}{2-p}}C_{*}^{\frac{2p}{p-2}}
+Q_{1}\right ]\\
&\quad\frac{ t^{\alpha }}{\alpha \Gamma(\alpha )}\\ &\leqslant
2\delta\left \| u_{0} \right \|_{L^{\infty }(\mathbb{R}^{N})}^{2}+\left [
2\delta\lambda^{*}\left \| u_{0} \right \|_{L^{\infty
}(\mathbb{R}^{N})}^{2}+\varepsilon^{\frac{2}{2-p}}C_{*}^{\frac{2p}{p-2}}
+Q_{1}\right ]\\
&\quad\frac{ T^{\alpha}}{\alpha \Gamma(\alpha )} \\ &=:M_{1},
 \end{aligned}
 $$
where $\lambda^{*}=(C_{2}-2\gamma )+\frac{p}{(2-p)\varepsilon
^{\frac{2}{p}}},\varepsilon >0$.\\
 \textbf{Case 2.}$N=2$. From \eqref{26}, we get
 \begin{equation}\label{2.10}
   \begin{aligned}
   &2\mu \int _{B(x,\delta )}u^{3}dy\leq 2\int _{B(x,\delta )}\left | \nabla u \right |^{2}dy+2\mu (1+\mu C_{GN}^{2}(2,\delta))\\
   &\qquad (\int _{B(x,\delta )}u^{3}dy\int _{B(x,\delta )}udy)+2\mu C_{GN}^{4}(2,\delta).
   \end{aligned}
 \end{equation}
Plugging \eqref{2.10} into \eqref{11} for $k\geq k^{*}=\frac{\mu
C_{GN}^{2}(2,\delta)+1}{\eta }$, we get
 \begin{align*}
    &\frac{\partial^{\alpha} }{\partial t^{\alpha}}\int _{B(x,\delta )}u^{2}dy++C_{*}^{-p}\left \| u \right \|_{L^{2}(B(x,\delta ))}^{p}+2\gamma\int _{B(x,\delta )}u^{2}dy\\
    &\leq C_{2}\int _{B(x,\delta )}u^{2}dy+2\mu C_{GN}^{4}(2,\delta).
\end{align*}
 Denote $w(t)=\int_{B(x,\delta )}u^{2}dy$ the solution of the following fractional differential equation
 \begin{equation}\label{2.12}
   \begin{cases}
     _{0}^{C}\textrm{D}_{t}^{\alpha }w(t)+(2\gamma -C_{2})w(t)+C_{*}^{-p}w^{\frac{p}{2}}(t)=2\mu C_{GN}^{4}(2,\delta),\\
   w(0)=(2\delta)^{2}\left \| u_{0} \right \|_{L^{\infty}(\mathbb{R}^{N})}^{2},
   \end{cases}
 \end{equation}
for any $(x,t)\in \mathbb{R}^{2}\times [0,T_{max})$. Then, we obtain
  $$
 \begin{aligned}
 &\int _{B(x,\delta )}u^{2}dy=\left \| u \right \|_{L^{2}(B(x,\delta ))}^{2}=  w(t)\\
&\leqslant w(0)+\left [
\lambda^{*}w(0)+\varepsilon^{\frac{2}{2-p}}C_{*}^{\frac{2p}{p-2}} +2\mu
C_{GN}^{4}(2,\delta)\right ]\\
&\quad\int_{0}^{t}(t-s)^{\alpha-1}E_{\alpha ,\alpha
}(\lambda^{*}(t-s)^{\alpha})ds\\ 
&\leqslant (2\delta)^{2}\left \| u_{0} \right
\|_{L^{\infty }(\mathbb{R}^{N})}^{2}\\
&+\left [ (2\delta)^{2}\lambda^{*}\left \|
u_{0} \right \|_{L^{\infty
}(\mathbb{R}^{N})}^{2}+\varepsilon^{\frac{2}{2-p}}C_{*}^{\frac{2p}{p-2}} +2\mu
C_{GN}^{4}(2,\delta)\right ]\frac{ T^{\alpha}}{\alpha \Gamma(\alpha )}\\
&=:M_{2}
 \end{aligned}
 $$
In summary, for any  $(x,t)\in \mathbb{R}^{N}\times [0,T_{max}]$, we have
 \begin{equation}\label{2.13}
   \left \| u \right \|_{L^{2}(B(x,\delta ))}\leq M=\begin{cases}
	\sqrt{M_{1}}\qquad N=1,\\
	\sqrt{M_{2}}\qquad N=2.
\end{cases}
 \end{equation}
 Then, we obtain
 \begin{equation}\label{2.14}
   \left \| u \right \|_{L^{1}(B(x,\delta ))}\leq (2\delta)^{\frac{N}{2}}M.
 \end{equation}
In order to  improve the $L^{2}$ boundedness of $u$ to $L^{\infty}$, by using
fractional Duhamel formula to equation \eqref{1}-\eqref{2} and let $f(u)=\mu
u^{2}(1-kJ*u)-\gamma u$, for all $(x,t)\in \mathbb{R}^{N}\times [0,T)$,
we have the solution
 $$u(t)=\mathcal{S}_{\alpha }(t)u_{0}+\int_{0}^{t}(t-s)^{\alpha -1}\mathcal{K} _{\alpha }(t-s)f(u(s))ds.$$
From Lemma \ref{lem3}-\ref{lemma2.6} and Proposition \ref{proposition.1}, we
have 
\begin{equation} 
\begin{aligned}
  \nonumber&0\leq u(x,t)\leq \left \| \mathcal{S} _{\alpha }(t)u_{0} \right \|_{L^{\infty }(\mathbb{R}^{N})}\\
  &\quad+\int_{0}^{t}(t-s)^{\alpha -1}\left \| \mathcal{K}_{\alpha }(t-s)f(u(s)) \right \|_{L^{\infty }(\mathbb{R}^{N})}ds\\
\nonumber&\leq C_{4}\left \| u_{0} \right \|_{L^{\infty
}(\mathbb{R}^{N})}+C_{4}\int_{0}^{t}(t-s)^{\alpha -1}\left \| f(u(s)) \right
\|_{L^{\infty }(\mathbb{R}^{N})}ds\\ \nonumber&\leq C_{4}\left \| u_{0} \right
\|_{L^{\infty }(\mathbb{R}^{N})}+\mu C_{4}\int_{0}^{t}(t-s)^{\alpha -1}\left \|
u^{2}(s) \right \|_{L^{\infty }(\mathbb{R}^{N})}ds\\ \nonumber&\leq  C_{4}\left
\| u_{0} \right \|_{L^{\infty }(\mathbb{R}^{N})}+\mu
C_{4}\int_{0}^{t}(t-s)^{\alpha -1}\left \| u(s) \right \|^{2}_{L^{\infty
}(\mathbb{R}^{N})}ds\\ \nonumber&\leq C_{4}\left \| u_{0} \right \|_{L^{\infty
}(\mathbb{R}^{N})}+\mu C_{4}M^{2}\int_{0}^{t}(t-s)^{\alpha -1}ds\\
\nonumber&\leq C_{4}\left \| u_{0} \right \|_{L^{\infty }(\mathbb{R}^{N})}+\mu
C_{4}M^{2}T^{\alpha}\frac{1}{\alpha}. \end{aligned} \end{equation} So there is
\begin{equation}\label{2.22} 0\leq u(x,t)\leq \left \| u_{0} \right
\|_{L^{\infty}(\mathbb{R}^{N})}+ \mu M^{2}C_{2}(N,T,\alpha), \end{equation}
with $M$ defined in \eqref{2.13} and where $$Q_{1}=2\mu \left ( (\mu ^{\frac{1}{3}}C_{GN}^{\frac{4}{3}}(1,\delta)+1)^{6}(\eta k)^{-5}+C_{GN}^{10}(1,\delta) \right )$$
  and $$\lambda^{*}=(C_{2}-2\gamma )+\frac{p}{(2-p)\varepsilon ^{\frac{2}{p}}},\varepsilon >0.$$
From Lemma \ref{Proposition2.1}, we know $T_{max}=+\infty$. As a conclusion, we
have shown that understanding $u$ in time is the existence of  weak solution in
Theorem \ref{theorem1}, and the blow-up criterion in Lemma \ref{Proposition2.1}
shows that $u$ is the unique weak solution for \eqref{1}-\eqref {2} on
$(x,t)\in\mathbb{R}^{N}\times [0,T)$ . Thus, the global boundedness of $u$ is proved.
\end{pro}
\begin{remark}
    Firstly, on the basis of the existence of the solution, Literature \cite{0Global} uses the principle of comparison to obtain the expression of the solution, and this paper uses the fractional  Duhamel formula (see Appendix \ref{90opliuyhhn}) to construct the solution. Secondly, in the process of proof, reference \cite{0Global} mainly uses local energy estimation (Gagliardo-Nirenberg inequality) to reconcile thermonuclear decomposition and so on to prove the bounded nature of the solution. In addition to local energy estimation (Gagliardo-Nirenberg inequality), this paper also uses fractional  differential equation, the symmetry of the kernel of 
    $(-\Delta)_{p}^{s}$ and other properties to verify. Finally, compared with reference \cite{0Global}, the blow-up for the solution is also presented in Appendix \ref{0po9754}.
\end{remark}

\section{Asymptotic behavior of weak solutions}\label{0kjhgfrt} 
Now, we consider the asymptotic behavior of the weak solution of  \eqref{1}-\eqref{2}.

 \begin{proposition}\label{proposition3.1}
Based on Theorem \ref{theorem4}, we have that $$\left \| u(x,t) \right
\|_{L^{\infty }(\mathbb{R}^{N}\times (0,T))}< a$$ holds, denoting $A,a$ in
Appendix \ref{i9ojjg}. Assume that the function
   $$H(x,t)=\int _{B(x,\delta)}h(u(y,t))dy$$
with
   $$h(u)=Aln\left ( 1-\frac{u}{A} \right )-aln\left ( 1-\frac{u}{a} \right )$$
is nonnegative and satisfies
   \begin{equation}\label{3.1}
     \frac{\partial^{\alpha} H(x,t)}{\partial t^{\alpha}}\leq -(-\Delta)_{p}^{s} H(x,t)+\int _{B(x,\delta)}\left | \nabla u(y,t) \right |^{2}dy-D(x,t)
   \end{equation}
with
   $$D(x,t)=\frac{1}{2}(A-a)\mu k\int _{B(x,\delta)}u^{2}(y,t)dy.$$
 \end{proposition}

 \begin{pro}
 Let $K=\left \| u(x,t) \right \|_{L^{\infty }(\mathbb{R}^{N}\times [0,T))}$, then $0<K<a$. From the definition of $h(\cdot)$, it is easy to deduce that
 \begin{equation}\label{3.2}
   {h}'(u)=\frac{a}{a-u}-\frac{A}{A-u}=\frac{(A-a)u}{(A-u)(a-u)},
 \end{equation}
and
 \begin{equation}\label{3.3}
  {h}''(u)=\frac{a}{(a-u)^{2}}-\frac{A}{(A-a)^{2}}=\frac{(Aa-u^{2})(A-a)}{(A-u)^{2}(a-u)^{2}}.
 \end{equation}
 By Lagrange mean value theorem and let $\left | {h}'(u(\xi )) \right |^{p-1}<1$, we will multiply \eqref{1} by ${h}'(u)\varphi _{\varepsilon }$ and integrate by parts over $B(x,\delta)$, where $\varphi _{\varepsilon }(\cdot )\in C_{0}^{\infty }(B(x,\delta)),\varphi _{\varepsilon }(\cdot )\rightarrow 1$ in $B(x,\delta)$ as $\varepsilon \rightarrow 0$. Then, we obtain
 \begin{align*}
     &(-\Delta )_{p}^{s}h(u)\\
     &=\int _{B(x,\delta)}\frac{\left | h(u(x))- h(u(y))\right |^{p-2}(h(u(x))-h(u(y)))}{\left | u(x)-u(y) \right |^{N+sp}}dy\\
&=\int _{B(x,\delta)}\frac{\left | {h}'(u(\xi )) \right |^{p-2}\left |
u(x)- u(y)\right |^{p-2}(u(x)-u(y)){h}'(u(\xi ))}{\left | x-y \right
|^{N+sp}\left | {u}'(\eta ) \right |^{N+sp}}dy\\ &\leq \left | {h}'(u(\xi ))
\right |^{p-1}{h}'(u)\left ( -\Delta  \right )_{p}^{s}u\leq {h}'(u)\left (
-\Delta  \right )_{p}^{s}u,
 \end{align*}
where $\xi ,\eta \in B(x,\delta )$, and 
\begin{align*}
    &(-\Delta )_{p}^{s }\int _{B(x,\delta)}h(u)dy\\
    &=\int _{B(x,\delta)}\frac{\left |\int _{B(x,\delta)}h(u(x))dx-\int _{B(x,\delta)}h(u(y))dx\right |^{p-2}}{\left | u(x)-u(y) \right |^{N+2s }}\\
    &\quad\int _{B(x,\delta )}h(u(x))-h(u(y))dxdy\\
&\leq \int _{B(x,\delta)}\int _{B(x,\delta)}\frac{\left
|h(u(x))-h(u(y))\right |^{p-2}(h(u(x))-h(u(y)))}{\left | u(x)-u(y) \right
|^{N+2s }}dxdy\\
&\leq \int _{B(x,\delta)}(-\Delta )^{s }h(u)dy. \end{align*}
From  equation \eqref{4}, we get
 $$\int _{B(x,\delta)}{h}'(u)\partial _{t}^{\alpha }udy\geq \int _{B(x,\delta)}\partial _{t}^{\alpha }h(u)dy=\partial _{t}^{\alpha }\int _{B(x,\delta)}h(u)dy.$$
Then
 \begin{align*}
& \frac{\partial^{\alpha} }{\partial t^{\alpha}}\int _{B(x,\delta)}h(u)\varphi
_{\varepsilon }dy\\
&\leq-\int _{B(x,\delta)}{h}'(u)(-\Delta)_{p}^{s} u\varphi
_{\varepsilon }dy\\
&\qquad+\int _{B(x,\delta)}[\mu u^{2}(1-kJ*u)-\gamma u]{h}'(u)\varphi
_{\varepsilon }dy\\
&\leq -(-\Delta )_{p}^{s}\int _{B(x,\delta)}\varphi
_{\varepsilon }h(u)dy\\
&\qquad+\int _{B(x,\delta)}[\mu u^{2}(1-kJ*u)-\gamma
u]{h}'(u)\varphi _{\varepsilon }dy.
 \end{align*}
 Taking $\varepsilon \rightarrow 0$, we obtain
 \begin{align*}
& \frac{\partial^{\alpha} }{\partial t^{\alpha}}\int _{B(x,\delta)}h(u)dy\\
&\leq-\int _{B(x,\delta)}(-\Delta)_{p}^{s} h(u)dy\\
&\qquad+\int _{B(x,\delta )}[\mu
u^{2}(1-kJ*u)-\gamma u]{h}'(u)dy\\ 
&\leq -(-\Delta )_{p}^{s }\int
_{B(x,\delta)}h(u)dy\\
&\qquad+\int _{B(x,\delta)}[\mu u^{2}(1-kJ*u)-\gamma u]{h}'(u)dy.
 \end{align*}
 which is
 \begin{align}\label{3.4}
\nonumber &\frac{\partial^{\alpha} }{\partial t^{\alpha}}H(x,t)\leq -(-\Delta )_{p}^{s }H(x,t)\\
 &\quad+\int _{B(x,\delta)}[\mu u^{2}(1-kJ*u)-\gamma u]{h}'(u)dy.
 \end{align}
By \eqref{1.7} in Appendix \ref{i9ojjg}, combining \eqref{3.5}-\eqref{3.7}, we
obtain 
\begin{equation}\label{3.8} \begin{aligned}
  &\int _{B(x,\delta)}{h}'(u)[\mu u^{2}(1-kJ*u)-\gamma u]dy\\
  &\leq -\frac{1}{2}(A-a)\mu k\int _{B(x,\delta)}u^{2}(y,t)dy\\
&\qquad+\frac{(A-a)k^{4}\mu k(2\delta )^{2}}{2(A-a)^{2}(a-k)^{2}}\int
_{B(x,\delta)}\left | \nabla u(y,t) \right |^{2}dy. \end{aligned}
\end{equation} 
From \eqref{3.3} and we obtain$${h}''(u)>
\frac{(A-a)^{2}}{A^{2}a},$$ inserting \eqref{3.8} into \eqref{3.4}, we obtain
\begin{equation}\label{3.9}
  \begin{aligned}
 & \frac{\partial^{\alpha} }{\partial t^{\alpha}}H(x,t)\leq -(-\Delta)_{p}^{s} H(x,t)+\frac{(A-a)k^{4}\mu k(2\delta )^{2}}{2(A-k)^{2}(a-k)^{2}}\\
  &\quad \int _{B(x,\delta)}\left | \nabla u(y,t) \right |^{2}dy-\frac{1}{2}(A-a)\mu k\int _{B(x,\delta)}u^{2}(y,t)dy.
  \end{aligned}
\end{equation} By choosing a sufficiently small $\delta$ such that
$$\frac{(A-a)k^{4}\mu k(2\delta )^{2}}{2(A-k)^{2}(a-k)^{2}}\leq 1.$$ Then, we
obtain $$\frac{\partial^{\alpha} }{\partial t^{\alpha}}H(x,t)\leq
-(-\Delta)_{p}^{s} H(x,t)+\int _{B(x,\delta)}\left | \nabla u(y,t) \right
|^{2}dy-D(x,t),$$ with $$D(x,t)=\frac{1}{2}(A-a)\mu k\int
_{B(x,\delta)}u^{2}(y,t)dy.$$
\end{pro}
\begin{remark}
    This proposition is based on Theorem \ref{theorem4}, comparing with Literature \cite{0Global}, we constructed the Lyapunov type function using the related properties of fractional $p$-Laplacian and using \eqref{3.5}-\eqref{3.7} in Appendix \ref{i9ojjg}.
\end{remark}

\begin{thm}\label{theorem5}
 Denote $u(x,t)$ the globally bounded solution of \eqref{1}-\eqref{2}.\\
 \textbf{(i).}\quad For any $\gamma\geq 1$, there exist $\mu
^{*}>0$ and $m^{*}>0$ such that for $\mu \in (0,\mu ^{*})$ and $\left \| u_{0}
\right \|_{L^{\infty }(\mathbb{R}^{N})}< m^{*}$, we have $$\left \| u(x,t)
\right \|_{L^{\infty }(\mathbb{R}^{N}\times [0,T))}< \frac{\gamma }{\mu
},$$ and therefore $$\left \| u(x,t) \right \|_{L^{\infty
}(\mathbb{R}^{N})}\leq \left \| u_{0} \right \|_{L^{\infty
}(\mathbb{R}^{N})}e^{(-\sigma)^{\frac{1}{\alpha }} t}.$$ where $\sigma :=\gamma
-\mu \left \| u(x,t) \right \|_{L^{\infty }(\mathbb{R}^{N}\times
[0,T))}>0.$

\textbf{(ii).}\quad  if $1\leq \gamma < \frac{\mu }{4k} $, there exist $\mu ^{**}>0$ and
$m^{**}>0$ such that for $\mu \in (0,\mu ^{**})$ and  $\left \| u_{0} \right
\|_{L^{\infty }(\mathbb{R}^{N})}< m^{**}$, we obtain  $$\left \| u(x,t) \right
\|_{L^{\infty }(\mathbb{R}^{N}\times [0,T))}< a,$$ and therefore
$$\lim_{t\rightarrow \infty }u(x,t)=0$$ locally uniformly in $\mathbb{R}^{N}$.

\end{thm}

\begin{pro}

 From  Theorem \ref{theorem4}, for any ${K}'>0$ and the definition
of $M$, there exist  $\mu ^{*}({K}')>0$ and $m_{0}({K}')>0$ such that for $\mu
\in  (0,\mu ^{*}({K}'))$ and $\left \| u_{0} \right \|_{L^{\infty
}(\mathbb{R}^{N})}<m_{0}({K}')$, there is $M$ which is sufficiently small such
that
 $$\left \| u(x,t) \right \|_{L^{\infty }(\mathbb{R}^{N}\times [0,T))}< {K}'.$$

\textbf{(i) The case:} $\left \| u(x,t) \right \|_{L^{\infty }(\mathbb{R}^{N}\times [0,T))}<{K}'=\frac{\gamma }{\mu }$.

Noticing $\sigma_{1}(t) =\gamma -\mu \left \| v(x,t) \right \|_{L^{\infty
}(\mathbb{R}^{N}\times [0,T))}>0$, let us consider the function
$v(x,t):=v(t)>0$ for all $x\in \mathbb{R}^{N}$. Then it follows that
 \begin{equation}\label{3.10}
   \begin{aligned}
   &\frac{\partial^{\alpha} v}{\partial t^{\alpha}}=-(-\Delta)_{p}^{s} v+\mu v^{2}(1-kJ*v)-\gamma v\\
&=-(-\Delta)_{p}^{s} v+v[\mu v(1-kJ*v)-\gamma ]\\ &=-(-\Delta)_{p}^{s}
v-v[\gamma -\mu v(1-kJ*v)]\\ &\leq -(-\Delta)_{p}^{s} v-v[\gamma -\mu v]\\
&\leq -(-\Delta)_{p}^{s} v-v\sigma_{1}.
   \end{aligned}
 \end{equation}
From Lemma \ref{Gao}, we know that $(-\Delta)_{p}^{s}v(t)\geq 0$. The last
expression can be rewritten in the following form
 \begin{equation}\label{dhjdhj}
   \begin{cases}
   _{0}^{C}\textrm{D}_{t}^{\alpha }v(t)=-\sigma_{1} v(t),\\
  v(0)=\left \| v_{0} \right \|_{L^{\infty }(\mathbb{R}^{N})},
   \end{cases}
 \end{equation}
 which ensures that $v(t)$ satisfies \eqref{1}-\eqref{2} with initial data $0< \underset{x\in \mathbb{R}^{N}}{max}\,u_{0}(x)\leq v_{0}$. From \eqref{9i9i} and \eqref{frfr},  we obtain that the solution to equation \eqref{dhjdhj} satisfies estimates
 $$v(t)\leq v_{0}E_{\alpha}(-\sigma t^{\alpha})\leq Cv_{0}e^{(-\sigma_{1})^{\frac{1}{\alpha }t} },\quad t>0$$
As $0<u_{0}\leq v_{0}$, then $v(t)$ is a supersolution of problem
\eqref{1}-\eqref{2}. Thus, by Lemma \ref{Gao}, we obtain
 $$u(x,t)\leq C \left \| u_{0} \right \|_{L^{\infty }(\mathbb{R}^{N})}e^{(-\sigma)^{\frac{1}{\alpha }} t}.$$
 For any $(x,t)\in \mathbb{R}^{N}\times [0,T)$, we get that
$$\left \| u(\cdot,t) \right \|_{L^{\infty }(\mathbb{R}^{N})}\leq C\left \|
u_{0} \right \|_{L^{\infty }(\mathbb{R}^{N})}e^{(-\sigma)^{\frac{1}{\alpha }}
t}.$$

  \textbf{(ii) The case:}\, $\left \| u(x,t) \right \|_{L^{\infty }(\mathbb{R}^{N}\times [0,T))}< {K}'=a$.
 
 First,  from \eqref{3.1} in Proposition \ref{proposition3.1}, we let $H_{0}(y)=H(y,0)$  and use fractional Duhamel formula, for any $(x,t)\in \mathbb{R}^{N}\times [0,T )$, then we get
 \begin{align*}
    & H(x,t)\leq \left \| H_{0} \right \|_{L^{\infty
}(\mathbb{R}^{N})}\\
&\quad+\int_{0}^{t}(t-s)^{\alpha -1}\mathcal{K}_{\alpha
}(t-s)(\left \| \nabla u \right \|_{2}^{2}-D(y,s))ds,
 \end{align*}
 from which we get
 \begin{align*}
     &\int_{0}^{t}(t-s)^{\alpha -1}\mathcal{K}_{\alpha }(t-s)D(y,s)ds \\
     &\leq \left \|
H_{0} \right \|_{L^{\infty }(\mathbb{R}^{N})}+\int_{0}^{t}(t-s)^{\alpha
-1}\mathcal{K}_{\alpha }(t-s)\left \| \nabla u \right \|_{2}^{2}ds.
 \end{align*}
 Due to
the definition of weak solution $u$, we obtain that $$\mathcal{K}_{\alpha
}(t-s)D(y,s)\in C^{2,1}(\mathbb{R}^{N}\times [0,T)),$$ This means that for all $x \in
\mathbb{R}^{N}$, the following restrictions apply: $$\lim_{t\rightarrow \infty
}\lim_{s\rightarrow t}\mathcal{K}_{\alpha }(t-s)D(y,s)=0,$$ or equivalently
$$\lim_{t\rightarrow \infty }\lim_{s\rightarrow t}\mathcal{K}_{\alpha
}(t-s)\int _{B(x,\delta)}u^{2}(z,s)dz=0,$$ Together with the fact that the heat
kernel converges to a delta function at $r\rightarrow t$, we can conclude that
for any $x\in \mathbb{R}^{N}$. $$\lim_{t\rightarrow \infty }\int
_{B(x,\delta)}u^{2}(y,t)dy=0.$$

 Moreover, due to the uniform boundedness of $u$ on $\mathbb{R}^{N}\times[0,T)$, by Lemma\ref{lemma5}, we can obtain the global boundedness of $\left \| \nabla u(\cdot ,t) \right \|_{L^{\infty }(\mathbb{R}^{N})}$. From this last point and Lemma \ref{lemma2.1}, we let $\Omega =B(x,\delta),p=q=\infty ,r=2$ and  the convergence of $\left \|  u(\cdot ,t)\right \|_{L^{\infty }(B(x,\delta)}$. That is, we get
$$\left \|  u(\cdot ,t)\right \|_{L^{\infty }(B(x,\delta))}\rightarrow 0$$ as
$t \rightarrow 0$. Thus, for any compact set in $\mathbb{R}^{N}$, by finite
coverage, we get that $u$ converges uniformly to $0$ in that compact set, which
means that when $t\rightarrow\infty$, we have that $u$ converges locally and
uniformly to $0$ in $\mathbb{R}^{N}$.
 The proof is complete.
 \end{pro}
\begin{remark}
    First, a local Lyapunov-type functional (Proposition \ref{proposition3.1}) will be constructed  before the proof on the sense of fractional derivative. And this paper verifies that \eqref{3.1} in Proposition \ref{proposition3.1} uses the fractional  Duhamel formula again. Secondly, we validate the comparison principle (see Lemma \ref{ollkkkmh}) of equation \eqref{1}-\eqref{2} more than in literature \cite{0Global}. Finally, based on the asymptotic convergence of the thermonuclear, the solution to the differential equation is changed to the fractional differential equation.
\end{remark}

 \section{Global boundness of solutions for \eqref{114}-\eqref{1.1.5}}\label{9o0oiyhnn}
 In the section, we will make $J=k=\mu=\gamma =1$ and replace fractional diffusion $(-\Delta)_{p}^{s}u$ as nonlinear fractional diffusion $(-\Delta)_{p}^{s}u^{m}$ in \eqref{1}-\eqref{2}  to get equation \eqref{114}-\eqref{1.1.5}. And we make $f(x,t)=u^{2}(1-\int _{\mathbb{R}^{N}}udx)-u,$ then equation \eqref{114}-\eqref{1.1.5} can write the following form
 \begin{align}\label{100}
\left\{\begin{matrix} \frac{\partial^{\alpha}u}{\partial t^{\alpha}} +\left (
-\Delta  \right )_{p}^{s}(u^{m})=f(x,t),& x\in \mathbb{R}^{N}\times(0,T] \\
 u(x,0)=u_{0}(x),& x\in \mathbb{R}^{N}.
 \end{matrix}\right.
\end{align} We take the announced values of the parameters $s,p$ with the
condition $sp<1$. Let $u\geq v$ be two ordered solutions of the equation
\eqref{100} defined in a strip $Q_{T}=\mathbb{R}^{N}\times(0,T)$. 

\begin{thm}\label{theorem6}
 \begin{align}
  &\frac{\partial^{\alpha } u}{\partial t^{\alpha }}+(-\Delta)_{p}^{s} u^{m}=u^{2}(1-\int _{\mathbb{R}^{N}}udx)-u \label{114}\\
&u(x,0)=u_{0}(x) \qquad x \in \mathbb{R}^{N}\label{1.1.5}
 \end{align}
with $  2-\frac{2}{N}<m\leq 3,N\geq 3,0<\alpha,s<1$ and $1<p<\frac{4}{3}$. Then
the  problem \eqref{114}-\eqref{1.1.5} has a unique weak solution which is
global bounded in fractional Sobolev space. then there exist $M>0,T>0$ such
that$$0\leq u(x,t)\leq M, \qquad (x,t)\in \mathbb{R}^{N}\times (0,T].$$
	\end{thm}
\begin{pro}

\textbf{Step 1.}\quad The $L^{r}$ estimates.

 For any $x\in \mathbb{R}^{N}$, multiply \eqref{114} by
$u^{k-1},k>1$ and integrating by parts over $\mathbb{R}^{N}$, by the proof of
Theorem 5.2 in \cite{2016fractional}, we obtain \begin{align*}
    &\int _{\mathbb{R}^{N}}u^{k-1}(-\Delta )_{p}^{s}u^{m}dx\\
&=\frac{1}{2}\int _{\mathbb{R}^{N}}\int
_{\mathbb{R}^{N}}\frac{(u^{m}(x,t)-u^{m}(y,t))\left | u^{m}(x,t)-u^{m}(y,t)
\right |^{p-2}}{\left | x-y \right |^{N+sp}}\\
&\qquad(u^{k-1}(x,t)-u^{k-1}(y,t))dxdy, \end{align*}
from Lemma \ref{233} and
equation \eqref{2.3}, let $a=u^{m}(x,t),b=u^{m}(y,t),\alpha =\frac{k-1}{m}$,
then we have 
\begin{align*}
    &\left | u^{m}(x,t)-u^{m}(y,t) \right |^{p-2}\left ( u^{m}(x,t) -u^{m}(y,t)\right )\\
   &\qquad \left ( u^{k-1}(x,t)-u^{k-1}(y,t) \right )\\
   &\geq c_{3}\left | u^{\frac{m(p-1)+k-1}{p}}(x,t)-u^{\frac{m(p-1)+k-1}{p}}(y,t)\right |^{p}.
\end{align*} 
By Sobolev inequality (Lemma \ref{1234}), we obtain that
\begin{align*}
   & \int _{\mathbb{R}^{N}}u^{k-1}(-\Delta )_{p}^{s}u^{m}dx\\
    &\geq  \frac{c_{3}}{2}\int _{\mathbb{R}^{N}}\int _{\mathbb{R}^{N}}\frac{\left | u^{\frac{m(p-1)+k-1}{p}}(x,t)-u^{\frac{m(p-1)+k-1}{p}}(y,t)\right |^{p}}{\left | x-y \right |^{N+sp}}dxdy\\
&\geq \frac{c_{3}S}{2}\left ( \int _{\mathbb{R}^{N}} \left |
u^{\frac{m(p-1)+k-1}{p}}(x,t)\right |^{p_{s }^{*}}dx\right
)^{\frac{p}{p_{s}^{*}}}, \end{align*} where $p_{s}^{*}=\frac{Np}{N-sp}$. From
Lemma \ref{lemma1}, we get $$\int
_{\mathbb{R}^{N}}u^{k-1}(_{0}^{C}\textrm{D}_{t}^{\alpha }u)dx\geq
\frac{1}{k}(_{0}^{C}\textrm{D}_{t}^{\alpha })\int _{\mathbb{R}^{N}}u^{k}dx.$$
Then, we obtain \begin{equation*}
  \begin{aligned}
&\frac{1}{k}(_{0}^{C}\textrm{D}_{t}^{\alpha }\int
_{\mathbb{R}^{N}}u^{k}dx)\leq\int
_{\mathbb{R}^{N}}u^{k-1}(_{0}^{C}\textrm{D}_{t}^{\alpha }u)dx\\ 
&=- \int
_{\mathbb{R}^{N}}u^{k-1}(-\Delta )_{p}^{s}u^{m}dx+\int
_{\mathbb{R}^{N}}u^{k+1}dx(1-\int _{\mathbb{R}^{N}}udx)\\
&\qquad-\int
_{\mathbb{R}^{N}}u^{k}dx\\ 
&\leq-\frac{c_{3}S}{2}\left ( \int _{\mathbb{R}^{N}}
\left | u^{\frac{m(p-1)+k-1}{p}}\right |^{p_{s }^{*}}dx\right
)^{\frac{p}{p_{s}^{*}}}\\
&\qquad+ \int _{\mathbb{R}^{N}}u^{k+1}dx(1-\int
_{\mathbb{R}^{N}}udx) -\int _{\mathbb{R}^{N}}u^{k}dx \end{aligned}
\end{equation*} 
thus, we have \begin{align}\label{4.1.1} \nonumber
&\frac{1}{k}(_{0}^{C}\textrm{D}_{t}^{\alpha }\int
_{\mathbb{R}^{N}}u^{k}dx)+\frac{c_{3}S}{2}\left ( \int _{\mathbb{R}^{N}} \left
| u^{\frac{m(p-1)+k-1}{p}}\right |^{p_{s }^{*}}dx\right
)^{\frac{p}{p_{s}^{*}}}\\
    &\leq\int _{\mathbb{R}^{N}}u^{k+1}dx(1-\int _{\mathbb{R}^{N}}udx)
-\int _{\mathbb{R}^{N}}u^{k}dx \end{align} 
In order to estimate $k\int
_{\mathbb{R}^{N}}u^{k+1}dx$, when $m\leq3$ and $k$ Satisfy \eqref{5.5}, then we
substitute \eqref{4.1.3}-\eqref{4.1.6} into \eqref{4.1.1} to get
\begin{align}\label{4.1.7} \nonumber &_{0}^{C}\textrm{D}_{t}^{\alpha }\int
_{\mathbb{R}^{N}}u^{k}dx+\frac{c_{3}Sk}{2}\left ( \int _{\mathbb{R}^{N}}
u^{\frac{m(p-1)+k-1}{p}p_{s}^{*}}dx\right )^{\frac{p}{p_{s}^{*}}}+k\int
_{\mathbb{R}^{N}}u^{k}dx\\ &\leq \frac{2k(k-1)}{(m+k-1)^{2}}\left \| \nabla
u^{\frac{k+m-1}{2}} \right \|_{L^{2}(\mathbb{R}^{N})}^{2}+C_{2}(N,k,m).
\end{align} 
Let $k\rightarrow \infty $, then $ \frac{2k(k-1)}{(m+k-1)^{2}}
\rightarrow 2$, so $\frac{2k(k-1)}{(m+k-1)^{2}} \leq 2.$  And from Lemma
\ref{5566} and $N\geq 3$, we obtain $$\frac{1}{S_{N}^{-1}\left \| u \right
\|_{1}^{\frac{1}{k-1}(2k/N+m-1)}}(\left \| u \right
\|_{k}^{k})^{1+\frac{m-1+2/N}{k-1}}\leq \left \| \nabla u^{(k+m-1)/2} \right
\|_{2}^{2}.$$ 
Let $k=\frac{N(1-m(p-1)}{sp}$, then
$\frac{m(p-1)+k-1}{p}p_{s}^{*}=k,\frac{p}{p_{s}^{*}}=\frac{k+m(p-1)-1}{k}$.
Therefore, we can get \begin{align*}
    &_{0}^{C}\textrm{D}_{t}^{\alpha }\int _{\mathbb{R}^{N}}u^{k}dx+\frac{c_{3}Sk}{2}\left ( \int _{\mathbb{R}^{N}}  u^{k}dx\right )^{\frac{k+m(p-1)-1}{k}}+k\int _{\mathbb{R}^{N}}u^{k}dx\\
&\leq \frac{2}{S_{N}^{-1}\left \| u \right
\|_{1}^{\frac{1}{k-1}(2k/N+m-1)}}(\left \| u \right
\|_{k}^{k})^{1+\frac{m-1+2/N}{k-1}}+C_{2}(N,k,m). \end{align*} Let $y(t)=\int
_{\mathbb{R}^{N}}u^{k}dx$ and $a=1+\frac{m-1+2/N}{k-1}$,
$$f(t)=\frac{2}{S_{N}^{-1}\left \| u \right
\|_{1}^{\frac{1}{k-1}(2k/N+m-1)}},\beta=\frac{k+m(p-1)-1}{k}.$$ Then, the
above-mentioned inequality can be written as \begin{equation}\label{2.27}
    _{0}^{C}\textrm{D}_{t}^{\alpha }y(t)+\frac{c_{3}Sk}{2}y^{\beta }(t)+ky(t)\leq f(t)y^{a}(t)+C_{2}(N,k,m)
\end{equation} By Lemma \ref{4554} and $0<\beta<a<1,t\in[0,T]$, fractional
differential inequality \eqref{2.27} has following solution
 \begin{align*}
   &y(t)\leq y(0)+(1-a)^{\frac{1}{1-\beta}}\varepsilon^{\frac{1}{1-a}}\frac{T^{\frac{\alpha}{1-\beta}}}{\alpha^{\frac{1}{1-\beta}}\Gamma^{\frac{1}{1-\beta}}(\alpha )}f^{\frac{1}{1-a}}(t)\\
     &\qquad+\left [\frac{[\lambda_{k} y^{1-\beta}(0)+(C_{2}(N,k,m)-\frac{c_{3}Sk}{2})(1-\beta)]T^{\alpha }}{\alpha \Gamma(\alpha )} \right ]^{\frac{1}{1-\beta}}.
\end{align*}
Therefore, we have
\begin{align}\label{4187} 
\nonumber  &y(t)=\int
_{\mathbb{R}^{N}}u^{k}dx \\
\nonumber&\leq y(0)+\left [\frac{[\lambda_{k}
y^{1-\beta}(0)+(C_{2}(N,k,m)-\frac{c_{3}Sk}{2})(1-\beta)]T^{\alpha }}{\alpha
\Gamma(\alpha )} \right ]^{\frac{1}{1-\beta}}\\
     &\qquad+(1-a)^{\frac{1}{1-\beta}}\varepsilon^{\frac{1}{1-a}}\frac{T^{\frac{\alpha}{1-\beta}}}{\alpha^{\frac{1}{1-\beta}}\Gamma^{\frac{1}{1-\beta}}(\alpha )}f^{\frac{1}{1-a}}(0).
\end{align} 
where $\lambda _{k}=-\frac{a-\beta }{\varepsilon ^{\frac{1-\beta
}{a-\beta }}}-k(1-\beta ),y(0)=\left \| u_{0} \right
\|_{L^{k}(\mathbb{R}^{N})}^{k}$ and $f(0)=\frac{2}{S_{N}^{-1}\left \| u_{0}
\right \|_{1}^{\frac{1}{k-1}(2k/N+m-1)}}$.

\textbf{Step 1.}\quad  The $L^{\infty}$ estimates.

 On account of the above arguments, our last task is to give the
uniform boundedness of solution for any $t>0$. Denote $q_{k}=2^{k}+2$, by
taking $k=q_{k}$ in \eqref{4.1.1}, we have \begin{equation}\label{4.2.1}
  \begin{aligned}
 &\frac{1}{q_{k}}(_{0}^{C}\textrm{D}_{t}^{\alpha }\int _{\mathbb{R}^{N}}u^{q_{k}}dx)+\frac{c_{3}S}{2}\left ( \int _{\mathbb{R}^{N}} \left | u^{\frac{q_{k}+m(p-1)-1}{p}}(x,t) \right |^{p_{s}^{*}}dx\right )^{\frac{p}{p_{s}^{*}}}\\
&\leq \int _{\mathbb{R}^{N}}u^{q_{k}+1}dx(1-\int _{\mathbb{R}^{N}}udx) -\int
_{\mathbb{R}^{N}}u^{q_{k}}dx.
  \end{aligned}
\end{equation}
By Lemma \ref{lemma333}, we substitute \eqref{4.2.2} into
\eqref{4.2.1} and with notice that $\frac{4q_{k}(q_{k}-1)}{(m+q_{k}-1)^{2}}\geq
2$. It follows \begin{equation}\label{4.2.3}
  \begin{aligned}
 &_{0}^{C}\textrm{D}_{t}^{\alpha }\int _{\mathbb{R}^{N}}u^{q_{k}}dx+\frac{c_{3}Sq_{k}}{2}\left ( \int _{\mathbb{R}^{N}} \left | u^{\frac{q_{k}+m(p-1)-1}{p}}(x,t) \right |^{p_{s}^{*}}dx\right )^{\frac{p}{p_{s }^{*}}}\\
 &\qquad+q_{k}\int _{\mathbb{R}^{N}}u^{q_{k}}dx\\
&\leq C(N)q_{k}^{\frac{\delta_{1}}{\delta_{1}-1}}(\int
_{\mathbb{R}^{N}}u^{q_{k-1}}dx)^{\gamma _{1}}+\frac{1}{2}\left \| \nabla
u^{\frac{m+q_{k}-1}{2}} \right \|_{2}^{2}\\ &\qquad+\frac{1}{2}\left \| u
\right \|_{L^{m+q_{k}-1}(\mathbb{R}^{N})}^{m+q_{k}-1}-\int
_{\mathbb{R}^{N}}udx\int _{\mathbb{R}^{N}}u^{q_{k}+1}dx
  \end{aligned}
\end{equation}
Invoking Lemma \ref{lemma333} once again, we sum up
\eqref{4.2.3} and \eqref{4.2.4}, with the fact that $\gamma _{1}\leq 2$ and
$\gamma _{2}<2$, we have 
\begin{equation*}
\begin{aligned}
&_{0}^{C}\textrm{D}_{t}^{\alpha }\int
_{\mathbb{R}^{N}}u^{q_{k}}dx+\frac{c_{3}Sq_{k}}{2}\left ( \int
_{\mathbb{R}^{N}} \left | u^{\frac{q_{k}+m(p-1)-1}{p}}(x,t) \right
|^{p_{s}^{*}}dx\right )^{\frac{p}{p_{s}^{*}}}\\
 &\qquad+q_{k}\int
_{\mathbb{R}^{N}}u^{q_{k}}dx\\ &\leq
C(N)q_{k}^{\frac{\delta_{1}}{\delta_{1}-1}}(\int
_{\mathbb{R}^{N}}u^{q_{k-1}}dx)^{\gamma _{1}}+\left \| \nabla
u^{\frac{m+q_{k}-1}{2}} \right \|_{L^{2}(\mathbb{R}^{N})}^{2}\\
 &\qquad+C_{3}(N)\\ &\leq
max\left \{ C(N),C_{3}(N)\right \}q_{k}^{\frac{\delta _{1}}{\delta
_{1}-1}}\\
 &\qquad\left [ (\int _{\mathbb{R}^{N}}u^{q_{k-1}}dx)^{\gamma 1}+1+\left \|
\nabla u^{\frac{m+q_{k}-1}{2}} \right \|_{L^{2}(\mathbb{R}^{N})}^{2}\right ] \\
&\leq 2max\left \{ C(N),C_{3}(N)\right \}q_{k}^{\frac{\delta _{1}}{\delta
_{1}-1}}\\
 &\qquad max\left \{ (\int _{\mathbb{R}^{N}}u^{q_{k-1}}dx)^{2},1,\left \| \nabla
u^{\frac{m+q_{k}-1}{2}} \right \|_{L^{2}(\mathbb{R}^{N})}^{2} \right \}.
\end{aligned} \end{equation*}
  Let $q_{k}=\frac{N(1-m(p-1)}{sp}$, then $\frac{m(p-1)+q_{k}-1}{p}p_{s}^{*}=q_{k},\frac{p}{p_{s}^{*}}=\frac{q_{k}+m(p-1)-1}{q_{k}}$. Therefore, we can get
\begin{align*}
    &_{0}^{C}\textrm{D}_{t}^{\alpha }\int _{\mathbb{R}^{N}}u^{q_{k}}dx+\frac{c_{3}Sq_{k}}{2}\left ( \int _{\mathbb{R}^{N}}  u^{q_{k}}dx\right )^{\frac{q_{k}+m(p-1)-1}{q_{k}}}\\
 &\qquad+q_{k}\int _{\mathbb{R}^{N}}u^{q_{k}}dx\\
    &\leq 2max\left \{ C(N),C_{3}(N)\right \}q_{k}^{\frac{\delta _{1}}{\delta _{1}-1}}\\
 &\qquad max\left \{ (\int _{\mathbb{R}^{N}}u^{q_{k-1}}dx)^{2},1,\left \| \nabla u^{\frac{m+q_{k}-1}{2}} \right \|_{L^{2}(\mathbb{R}^{N})}^{2} \right \}.
\end{align*}
   Let $$K_{0}=max\left \{ 1,\left \| u_{0} \right \|_{L^{1}(\mathbb{R}^{N})},\left \| u_{0} \right \| _{L^{\infty }(\mathbb{R}^{N})},\left \| \nabla u_{0}^{\frac{m+q_{k}-1}{2}} \right \|_{L^{2}(\mathbb{R}^{N})}^{2}\right \},$$
   we have the following inequality for initial data
\begin{align*}
  &\int _{\mathbb{R}^{N}}u_{0}^{q_{k}}dx\\
 &\leq \left ( max\left \{ \left \| u_{0} \right \|_{L^{1}(\mathbb{R}^{N})},\left \| u_{0} \right \| _{L^{\infty }(\mathbb{R}^{N})},\left \| \nabla u_{0}^{\frac{m+q_{k}-1}{2}} \right \|_{L^{2}(\mathbb{R}^{N})}^{2}\right \}\right )^{q_{k}}\\
 &\leq K_{0}^{q_{k}}.
\end{align*}
  Let $d_{0}=\frac{\delta _{1}}{\delta _{1}-1}$, it is easy to that $q_{k}^{d_{0}}=(2^{k}+2)^{d_{0}}\leq (2^{k}+2^{k+1})^{d_{0}}$. By taking $\bar{a}=max\left \{ C(N),C_{3}(N) \right \}3^{d_{0}}$ in the Lemma \ref{lemma444}, we obtain
\begin{align}\label{4.2.6}
 \nonumber \int u^{q_{k}}dx&\leq (2\bar{a})^{2^{k}-1}2^{d_{0}(2^{k+1}-k-2)}\\
  &max\left \{ \underset{t\geq 0}{sup}(\int _{\mathbb{R}^{N}}u^{q}dx)^{2^{k}},k_{0}^{q_{k}} \right \}\frac{T^{\alpha }}{\alpha \Gamma (\alpha)}.
\end{align}
  Since $q_{k}=2^{k}+2$ and taking the power $\frac{1}{q_{k}}$ to both sides of \eqref{4.2.6}, then the boundedness of the solution $u(x,t)$ is obtained by passing to the limit $k\rightarrow\infty$
 \begin{equation*}
   \left \| u(x,t) \right \|_{L^{\infty }(\mathbb{R}^{N})}\leq 2\bar{a}2^{2d_{0}}max\left \{ \underset{t\geq 0}{sup}\int _{\mathbb{R}^{N}} u^{q_{0}}dx,K_{0}\right \}\frac{T^{\alpha }}{\alpha \Gamma (\alpha)}.
 \end{equation*}
  On the other hand, by \eqref{4187} with $q_{0}>2,t\in [0,T]$, we know
 \begin{align*}
  &\int _{\mathbb{R}^{N}}u^{q_{0}}dx\leq  \int _{\mathbb{R}^{N}}u^{3}dx\\
  &\leq \left [\frac{[\lambda_{k} y^{1-\beta}(0)+(C_{2}(N,k,m)-\frac{3c_{3}S}{2})(1-\beta)]T^{\alpha }}{\alpha \Gamma(\alpha )} \right ]^{\frac{1}{1-\beta}}\\
     &\qquad+y(0)+(1-a)^{\frac{1}{1-\beta}}\varepsilon^{\frac{1}{1-a}}\frac{T^{\frac{\alpha}{1-\beta}}}{\alpha^{\frac{1}{1-\beta}}\Gamma^{\frac{1}{1-\beta}}(\alpha )}f^{\frac{1}{1-a}}(0).
\end{align*} 
where $\lambda _{k}=-\frac{a-\beta }{\varepsilon ^{\frac{1-\beta
}{a-\beta }}}-3(1-\beta )$, and
$$y(0)=\left \| u_{0} \right
\|_{L^{3}(\mathbb{R}^{N})}^{3},f(0)=\frac{2}{S_{N}^{-1}\left \| u_{0} \right
\|_{1}^{\frac{1}{2}(6/N+m-1)}}.$$ 
Therefore we finally have
 \begin{align*}
   &\left \| u(x,t) \right \|_{L^{\infty }(\mathbb{R}^{N})}\\
   &\leq C(N,\left \| u_{0} \right \|_{L^{1}(\mathbb{R}^{N})},\left \| u_{0} \right \|_{L^{\infty }(\mathbb{R}^{N})},\left \| \nabla u_{0}^{\frac{m+2}{2}} \right \|_{L^{2}(\mathbb{R}^{N})}^{2},T^{\alpha})=M.
 \end{align*}
 The proof of this theorem is complete.
 \end{pro}
\begin{remark}
    Firstly, compared with Literature \cite{0Global}, it proves the global boundedness of solution in arbitrary dimension space, while this paper mainly proves in fractional  Sobolev space. Secondly, Literature \cite{0Global} mainly uses $L^{r}$ estimation, Moser iteration and Young inequality, Holder inequality, Gagliardo-Nirenberg inequality and other common inequalities to verify the equation. This paper is also based on the above-mentioned inequalities and the application of various fractional  differential inequalities in Appendix \ref{0opliutv}. Finally, the existence of the solution of equation \eqref{114}-\eqref{1.1.5} is discussed in Appendix \ref{28uijybb}.
\end{remark}

\section{Acknowledgements} 
This work is supported by the State Key Program of
National Natural Science of China under Grant No.91324201. This work is also
supported by the Fundamental Research Funds for the Central Universities of
China under Grant 2018IB017, Equipment Pre-Research Ministry of Education Joint
Fund Grant 6141A02033703 and the Natural Science Foundation of Hubei Province
of China under Grant 2014CFB865. 
\begin{appendices} 
\section{ Definitions, Related Lemma, and complements}

\subsection{Some useful Lemma about  fractional derivative } \label{0opliutv}
In the section, let us recall some necessary Lemma and useful properties of fractional
derivative. 
\begin{lemma}\textsuperscript{\cite{Kilbas2006}}\label{gugui}
  If $0<\alpha<1,u\in AC^{1}[0,T]$ or $u\in C^{1}[0,T]$,then the equality
  $$_{0}\textrm{I}_{t}^{\alpha }(_{0}^{C}\textrm{D}_{t}^{\alpha }u)(t)=u(t)-u(0)$$
  and
  $$_{0}^{C}\textrm{D}_{t}^{\alpha }(_{0}\textrm{I}_{t}^{\alpha }u)(t)=u(t),$$
hold almost everywhere on $[0,T]$, In addition,
\begin{align*}
   &_{0}^{C}\textrm{D}_{t}^{1-\alpha }\int_{0}^{t}\,_{0}^{C}\textrm{D}_{\tau
}^{\alpha }u(\tau )d\tau \\
&=\left ( _{0}\textrm{I}_{t}^{\alpha
}\frac{d}{dt}\,_{0}\textrm{I}_{t}^{1}\,_{0}\textrm{I}_{t}^{1-\alpha
}\frac{d}{dt}u \right )(t)=u(t)-u(0). 
\end{align*}
\end{lemma}

\begin{lemma}\textsuperscript{\cite{1948The}}\label{lemma288}
   If $0<\alpha<1,\eta>0$, then there is $0\leq E_{\alpha,\alpha }(-\eta )\leq \frac{1}{\Gamma (\alpha )}$. In addition, for $\eta>0$, $E_{\alpha ,\alpha }(-\eta )$ is a monotonically decreasing function.
\end{lemma}

\begin{lemma}\textsuperscript{\cite{Ahmed2017A}}\label{lemma23} Let
$0<\alpha<1$ and $v\in C([0,T],\mathbb{R}^{N}),{v}'\in
L^{1}(0,T;\mathbb{R}^{N})$ and $u$ be monotone. Then
  \begin{equation}\label{888}
    v(t)\partial _{t}^{\alpha }v(t)\geq \frac{1}{2}\partial _{t}^{\alpha }v^{2}(t),\quad t\in(0,T].
  \end{equation}
\end{lemma}

\begin{lemma}\textsuperscript{\cite{Ahmed2017A}}\label{lemma1} Let $0<\alpha<1$
and $u\in C([0,T],\mathbb{R}^{N}),{u}'\in L^{1}(0,T;\mathbb{R}^{N})$ and $u$ be
monotone. And when $n\geq 2$,the Caputo fractional derivative with respect to
time $t$ of $u$ is defined by \eqref{32}. Then there is \begin{equation*}
   u^{n-1}(_{0}^{C}\textrm{D}_{t}^{\alpha }u) \geq \frac{1}{n}(_{0}^{C}\textrm{D}_{t}^{\alpha }u^{n}).
\end{equation*} \end{lemma}

\begin{lemma}\textsuperscript{\cite{2016Weak}}\label{lemma2.14}
 Suppose that a nonnegative function $u(t)\geq 0$  satisfies
 \begin{equation}\label{4433}
   _{0}^{C}\textrm{D}_{t}^{\alpha }u(t)+c_{1}u(t)\leq f(t)
 \end{equation}
  for almost all $t\in [0,T]$, where $c_{1}>0$, and the function $f(t)$ is nonnegative and integrable for $t\in [0,T]$. Then
  \begin{equation}\label{55}
    u(t)\leq u(0)+\frac{1}{\Gamma (\alpha )}\int_{0}^{t}(t-s)^{\alpha -1}f(s)ds.
  \end{equation}
 \end{lemma}

 \begin{lemma}\textsuperscript{\cite{hui2022global}}\label{lemma444}
   Assume the function $y_{k}(t)$ is nonnegative and exists the Caputo fractional derivative for $t\in [0,T]$ satisfying
   \begin{align}\label{4.17}
  \nonumber _{0}^{C}\textrm{D}_{t}^{\alpha }y_{k}(t)&\leq -C_{9}(y_{k}(t))^{\frac{k+m-1}{k}}-y_{k}+a_{k}(y_{k-1}^{\gamma _{1}}(t)\\
  &\quad+y_{k-1}^{\gamma _{2}}(t)),
   \end{align}
   where $a_{k}=\bar{a}3^{rk}>1$ with $\bar{a},r$ are positive bounded constants and $0<\gamma _{2}<\gamma _{1}\leq 3$. Assume also that there exists a bounded constant $K\geq 1$ such that $y_{k}(0)\leq K^{3^{k}}$, then
   \begin{align}\label{4.18}
 \nonumber    y_{k}(t)&\leq (2\bar{a})^{\frac{3^{k}-1}{2}}3^{r (\frac{3^{k+1}}{4}-\frac{k}{2}-\frac{3}{4})}max\left \{ \underset{t\in [0,T]}{sup} y_{0}^{3^{k}}(t),K^{3^{k}}\right \}\\
    &\qquad \frac{T^{\alpha }}{\alpha \Gamma (\alpha)}.
   \end{align}
 \end{lemma}

\begin{lemma}\textsuperscript{\cite{hui2022global}}\label{4554}
 Suppose $0<k<m<1$ and $b(t)$ is continuous and boundary. And let $y(t)\geq 0$ be a solution of the fractional differential inequality
   \begin{equation}\label{256}
       _{0}^{C}\textrm{D}_{t }^{\alpha }y(t)+\alpha y^{k}(t)+\beta y(t)\leq b(t) y^{m}(t)+c_{4}.
   \end{equation}
   For almost all $t\in [0,T]$, then
   \begin{align*}
   y(t)&\leq y(0)+\left [\frac{[\lambda_{k} y^{1-k}(0)+(c_{4}-\alpha)(1-k)]T^{\alpha }}{\alpha \Gamma(\alpha )} \right ]^{\frac{1}{1-k}}\\
     &\qquad+(1-m)^{\frac{1}{1-k}}\varepsilon^{\frac{1}{1-m}}\frac{T^{\frac{\alpha}{1-k}}}{\alpha^{\frac{1}{1-k}}\Gamma^{\frac{1}{1-k}}(\alpha )}b^{\frac{1}{1-m}}(t).
\end{align*} where $\lambda _{k}=-\frac{m-k}{\varepsilon
^{\frac{1-k}{m-k}}}-\beta(1-k)$ and $\alpha,\beta,c_{4},\varepsilon>0$ are all
contants. \end{lemma}
\begin{lemma}\textsuperscript{\cite{2001Fractional}}\label{Gao}
  Let us consider the  fractional differential equation
  \begin{align}\label{jiojio}
  \begin{cases}
     _{0}^{C}\textrm{D}_{t}^{\alpha }u(t)=-wu(t),\,0<\alpha<1,w>0,\\
u(0)=u_{0}.
  \end{cases}
 \end{align}
Then, the solution of \eqref{jiojio} can be obtained by applying the Laplace
transform technique which implies: \begin{equation}\label{9i9i}
    u(t)=u_{0}E_{\alpha}(-wt^{\alpha}),\quad t>0.
\end{equation} \end{lemma} \begin{lemma}\textsuperscript{\cite{2001Fractional}}
   If $0<\alpha<1,t>0,w>0$, for Mittag-Leffler function $E_{\alpha,1 }(wt^{\alpha })$ , then there is a constant $C$ such that
   \begin{equation}\label{frfr}
     E_{\alpha}(wt^{\alpha})=E_{\alpha,1 }(wt^{\alpha })\leq Ce^{w^{\frac{1}{\alpha }t} }.
   \end{equation}
\end{lemma}

\begin{thm}\textsuperscript{\cite{2016fractional}}\label{1234} (Fractional
Sobolev inequality) Assume that $0<s<1$ and $p>1$ are such that $ps<N$. Then,
there exist a positive constant $S\equiv S(N,s,p)$ such that for all $v\in
C_{0}^{\infty }(\mathbb{R}^{N})$, 
\begin{align*}
    \left [ v \right ]_{W^{s,p}(\mathbb{R}^{N})}^{p}&=\iint_{\mathbb{R}^{2N}}\frac{\left | v(x)-v(y) \right |^{p}}{\left | x-y \right |^{N+sp}}dxdy\\
    &\geq S\left ( \int _{\mathbb{R}^{N}} \left | v(x) \right |^{p_{s}^{*}}\right )^{\frac{p}{p_{s}^{*}}},
\end{align*} 
where $p_{s}^{*}=\frac{pN}{N-sp}$. 
\end{thm}
\subsection{ Fractional Duhamel's formula }\label{90opliuyhhn}
According to \cite{giacomoni2018existence}, we define the operator $A:=-\left ( -\Delta
\right )_{p}^{s}$ on $L^{\infty}(\Omega )$ by \begin{equation}\label{scsc}
    \begin{cases}
      D(A)=\left \{ u\in X_{0}\cap L^{\infty }(\Omega ):Au\in L^{\infty }(\Omega ) \right \},\\
Au=-\left ( -\Delta  \right )_{p}^{s}u\quad \text{for}\quad u\in D(A).
    \end{cases}
\end{equation} It is well known that $A$ is m-accretive in $L^{\infty}(\Omega
)$. Then, we convert Equation \eqref{1}-\eqref{2} to the following abstract form:
\begin{equation}\label{fgfga} \begin{cases}
 \frac{\partial^{\alpha }u}{\partial t^{\alpha }}=Au+f(u),&(x,t)\in\Omega\times(0,T),\\
u(x,0)=u_{0}(x),& x\in\Omega, \end{cases} \end{equation} where $f(u(x,t))=\mu
u^{2}(1-kJ*u)-\gamma u$. Let $X_{1}$ be a Hilbert space and the function
$f:[0,T]\rightarrow X_{1}$ and $A$ is self-adjoint  on $X_{1}$. Then there
exists a measure space $(\sum ,\mu )$ and a Borel measurable function $a$ and a
unitary map $U:L^{\infty}(\sum ,\mu )\rightarrow X$ such that
$$U^{-1}AU=T_{a},\quad \text{where}\,T_{a}\varphi (\xi )=a(\xi )\varphi (\xi
),\quad\xi \in \sum $$

\begin{lemma}\textsuperscript{\cite{2021zhan}}\label{lem3}
   If $u$ satisfies Equation \eqref{fgfga} and $t\in [0,T]$, then $u$ also satisfies
   $$
   \begin{aligned}
   u(t)&=U(E_{\alpha }(a(\xi )t^{\alpha })U^{-1}u_{0})\\
   &\qquad+\int_{0}^{t}(t-s)^{\alpha -1}U(E_{\alpha ,\alpha }((t-s)^{\alpha }a(\xi ))U^{-1}f(u(s)))ds.
    \end{aligned}
   $$
 \end{lemma}
\begin{remark} By define two maps $$\mathcal{S} _{\alpha }(t)\phi=U(E_{\alpha
}(a(\xi )t^{\alpha })U^{-1}\phi)$$
and
$$ \mathcal{K} _{\alpha }(t)\phi=U(E_{\alpha
,\alpha }((t)^{\alpha }a(\xi ))U^{-1}\phi),\, \phi \in X,$$
then through Lemma
\ref{lem3}, we can get the following fractional Duhamel's formula. \end{remark}

\begin{lemma}\textsuperscript{\cite{2021zhan}}\label{lemma2.6} (fractional
Duhamel's formula)  If $u\in C([0,T];X)$ satisfies Equation \eqref{1}
-\eqref{2}, then $u$ satisfies the following integral equation:
   $$u(t)=\mathcal{S}_{\alpha }(t)u_{0}+\int_{0}^{t}(t-s)^{\alpha -1}\mathcal{K} _{\alpha }(t-s)f(u(s))ds.$$
   The solution operators $\mathcal{S}_{\alpha }(t)$ and $\mathcal{K} _{\alpha }(t)$ are defined by the functional calculus of $A$ via the Mittag-Leffler function when evaluated at $A$.
\end{lemma}
\begin{proposition}\textsuperscript{\cite{zhou2018duhamel}}\label{proposition.1}
   For each fixed $t \geq 0$,$\mathcal{S}_{\alpha }(t)$ and $\mathcal{K} _{\alpha }(t)$ are linear and bounded operators, for any $\phi \in X$,

   \begin{equation}\label{9}
     \left \|\mathcal{S}_{\alpha }(t)\phi  \right \|_{X}\leq C_{4}\left \| \phi \right \|_{X},\qquad  \left \|\mathcal{K}_{\alpha }(t)\phi  \right \|_{X}\leq C_{4}\left \| \phi \right \|_{X},
   \end{equation}
where $C_{4}$ is a constant.
 \end{proposition}

\begin{lemma}\textsuperscript{\cite{2017GagliardO}}\label{lemma99}
  Assume $1\leq p<N$, then $\forall f\in C_{c}^{\infty }(\mathbb{R}^{N})$, there is
  \begin{equation}\label{99}
      \left \| f \right \|_{L^{q}(\mathbb{R}^{N})}\leq C_{1}\left \| \nabla f \right \|_{L^{p}(\mathbb{R}^{N})}.
  \end{equation}
  Which is $q=\frac{Np}{N-p}$, and $C_{1}$ only dependent on $p,N$.
\end{lemma}
\subsection{Blow-up for equation \eqref{1}-\eqref{2}}\label{0po9754}
\begin{definition}\textsuperscript{\cite{2019the}}\label{4567}
  Suppose eigenfunction $e_{1}>0$ associated to the first eigenvalue $\lambda _{1}>0$ that satisfies the fractional eigenvalue problem
  \begin{align}\label{123456}
      &(-\Delta )_{p}^{s}e_{1}(x)=\lambda _{1}e_{1}(x),x\in \Omega ,\\
\nonumber &e_{1}(x)=0,\quad x\in \mathbb{R}^{N}\setminus  \Omega ,
  \end{align}
  normamlized such that $\int _{\Omega }e_{1}(x)dx=1$.
\end{definition}

\begin{lemma}\label{Proposition2.1} (Blow-up) Assume the initial data
$0<u_{0}\in X_{0}<1$.  And if
$1+\lambda _{1}\leq \int _{\Omega }u_{0}(x)e_{1}(x)dx=H_{0},$ and $\lambda
_{1},e_{1}(x)$ is the value in Definition \ref{4567}, then the solution of
problem \eqref{1}-\eqref{2} blow-up in a finite time
  $T_{max}$ that satisfies the bi-lateral estimate
  $$\left ( \frac{\Gamma (\alpha +1)}{4(H_{0}+1/2)} \right )^{\frac{1}{\alpha }}\leq T_{max}\leq \left ( \frac{\Gamma (\alpha +1)}{H_{0}} \right )^{\frac{1}{\alpha }}.$$
  Furthermore, if $T_{max}<\infty$, then
  $${\lim_{t\to T_{max}}}\Arrowvert u(\cdot,t)\Arrowvert_{L^{\infty}(\mathbb{R}^{N})}=\infty.$$
\end{lemma} 
\begin{pro} Multiplying equations \eqref{1} by $e_{1}(x)$ and
integrating over $\Omega $, we obtain 
\begin{align}
 \nonumber  & \partial _{t}^{\alpha }\int _{\Omega }u(x,t)e_{1}(x)dx+\int _{\Omega }(-\Delta )_{p}^{s}u(x,t)e_{1}(x)dx\\
&=\int _{\Omega }[\mu (1-kJ*u)u(x,t)-\gamma ]u(x,t)e_{1}(x)dx. \end{align} 
By \eqref{123456}, we have 
\begin{align*}
    \int _{\Omega }(-\Delta
)_{p}^{s}u(x,t)e_{1}(x)dx&=\int _{\Omega }u(x,t)(-\Delta
)_{p}^{s}e_{1}(x)dx\\
&=\lambda _{1}\int _{\Omega }u(x,t)e_{1}(x)dx,
\end{align*}
as
$u=0,e_{1}(x)=0,x\in \mathbb{R}^{N}\setminus  \Omega$ and $$\left ( \int
_{\Omega }u(x,t)e_{1}(x)dx \right )^{2}\leq \int _{\Omega
}u^{2}(x,t)e_{1}(x)dx,$$ let function $H(t)=\int _{\Omega }u(x,t)e_{1}(x)dx$,
then satisfies \begin{equation}
    \partial _{t}^{\alpha }H(t)+(\gamma +\lambda _{1})H(t)\geq \mu H^{2}(t).
\end{equation} The proof process later can refer to the reference
[\cite{2020Global11},Theorem 2.1] and [\cite{2014Blowing}, Theorem 3.2]. The
proof of the lemma is complete. 
\end{pro}
\begin{remark}
    The proof of the lemma is based on references \cite{2020Global11} and \cite{2014Blowing}, but the difference is that the properties of eigenfunction for fractional $p$-Laplacian and the correlation inequality are used to verify them.
\end{remark}

\section{Some useful results}\label{i9ojjg}

 To study the long time behavior of solutions for
\eqref{1}-\eqref{2}, by \cite{0Global}, we denote $$F(u):=\mu
u^{2}(1-kJ*u)-\gamma u.$$ For $1\leq\gamma<\frac{\mu }{4k}$, there are three
constant solutions for $F(u)=0$: $0,a,A$, where \begin{equation}\label{1.7}
  a=\frac{1-\sqrt{1-4k\frac{\gamma }{\mu }}}{2k},\qquad A=\frac{1+\sqrt{1-4k\frac{\gamma }{\mu }}}{2k},
\end{equation}
 and satisfy $1<\frac{\gamma }{\mu }<a<A$.
\begin{remark}
     By \eqref{1.7}, we can get $$\mu u^{2}(1-ku)-\gamma u=k\mu u(A-u)(u-a),$$
and 
\begin{equation}\label{3.5} 
\begin{aligned} &\int _{B(x,\delta)}{h}'(u)[\mu
u^{2}(1-kJ*u)-\gamma u]dy\\ &=\int _{B(x,\delta)}{h}'(u)[\mu u^{2}(1-ku)-\gamma
u]dy\\
&\quad+\mu k\int _{B(x,\delta)}{h}'(u)u^{2}(u-J*u)dy\\ &=-(A-a)\mu k\int
_{B(x,\delta)}u^{2}(y,t)dy\\
&\quad+\mu k\int _{B(x,\delta)}{h}'(u)u^{2}(u-J*u)dy.
\end{aligned} \end{equation} Noticing that when $0\leq u\leq k$, there is
$$0\leq {h}'(u)u\leq \frac{(A-a)k^{2}}{(A-k)(a-k)}.$$
From Young's inequality
and the median value theorem, we can get \begin{equation}\label{3.6}
\begin{aligned}
&\mu k\int _{B(x,\delta)}{h}'(u)u^{2}(u-J*u)dy\\
&\leq \mu
k\int _{B(x,\delta)}\int
_{B(x,\delta)}{h}'(u)u^{2}(y,t)(u(y,t)-u(z,t))J(y-z)dzdy\\ &\leq
\frac{(A-a)K^{2}}{(A-K)(a-K)}\mu k\int _{B(x,\delta)}\int
_{B(x,\delta)}u(y,t)\left | (u(y,t)-u(z,t)) \right |\\
&\qquad J(y,z)dzdy\\
&\leq\frac{(A-a)K^{4}\mu k}{2(A-K)^{2}(a-K)^{2}}\int _{B(x,\delta)}\int
_{B(x,\delta)}(u(z,t)-u(y,t))^{2}\\
&\qquad J(z-y)dzdy+\frac{1}{2}(A-a)\mu k\int
_{B(x,\delta)}u^{2}(y,t)dy\\ &\leq \frac{(A-a)K^{4}\mu
k}{2(A-K)^{2}(a-K)^{2}}\int _{B(x,\delta)}\int _{B(x,\delta)}\int_{0}^{1}\left
| \nabla u(y+\theta (z-y),t) \right |^{2}\\ &\qquad\left | z-y \right |^{2}
J(z-y) d\theta dzdy\\
&\qquad+\frac{1}{2}(A-a)\mu k\int _{B(x,\delta)}u^{2}(y,t)dy,
  \end{aligned}
\end{equation} changing the variables ${y}'=y+\theta (z-y) , {z}'=z-y$,
then$$\begin{vmatrix}
	& \frac{\partial{y}' }{\partial y} \qquad \frac{\partial {y}'}{\partial z}\\
	& \frac{\partial {z}'}{\partial y}\qquad  \frac{\partial {z}'}{\partial z}\\
\end{vmatrix}=\begin{vmatrix}
	1-\theta   \qquad\theta & \\
	-1\qquad	1 &
\end{vmatrix}=1-\theta+\theta=1.$$ For any $\theta \in  [0,1],y,z\in
B(x,\delta)$, we have ${y}'\in B((1-\theta )x+\theta z,(1-\theta )\delta )
,{z}'\in B(x-y,\delta )$. Noticing $B((1-\theta )x+\theta z,(1-\theta )\delta
)\subseteq B(x,\delta)$ and $B(x-y,\delta)\subseteq B(0,2\delta)$, we obtain
\begin{equation}\label{3.7} \begin{aligned} &\frac{(A-a)K^{4}\mu
k}{2(A-K)^{2}(a-K)^{2}}\int _{B(x,\delta)}\int _{B(x,\delta)}\int_{0}^{1}\left
| \nabla u(y+\theta (z-y),t) \right |^{2}\\ & \qquad\left | z-y \right |^{2}
J(z-y)d\theta dzdy\\ 
&\leq \frac{(A-a)K^{4}\mu
k}{2(A-K)^{2}(a-K)^{2}}\int_{0}^{1}d\theta \int _{B(x,\delta)}\int
_{B(x,\delta)}\left | \nabla u({y}',t) \right |^{2}\\ & \qquad\left | {z}' \right
|^{2}J({z}')d{z}'d{y}'\\
&\leq \frac{(A-a)K^{4}\mu k (2\delta
)^{2}}{2(A-K)^{2}(a-K)^{2}}\int _{B(x,\delta)}\left | \nabla u(y,t) \right
|^{2}dy. \end{aligned} \end{equation} \end{remark}
\begin{lemma}\textsuperscript{\cite{1966}}\label{lemma2.1}
  Let $\Omega$ be an open subset of $\mathbb{R}^{N}$, assume that $1\leq p,q \leq \infty$ with $(N-q)p<Nq $ and $r\in(0,p)$. Then there exists constant $C_{GN}>0$ only depending on $q, r$ and $\Omega$ such that for any $u\in W^{1.q}(\Omega)\cap L^{p}(\Omega)$
  $$\int_{\Omega}u^{p}dx\leq C_{GN}(\Arrowvert\nabla u\Arrowvert^{(\lambda^{*}p)}_{L^{q}(\Omega)}\Arrowvert u\Arrowvert_{L^{p}(\Omega)}^{(1-\lambda^{*})p}+\Arrowvert u\Arrowvert_{L^{p}(\Omega)}^{p})$$
  holds with $$\lambda^{*}=\frac{\frac{N}{r}-\frac{N}{p}}{1-\frac{N}{q}+\frac{N}{r}}\in (0,1).$$
\end{lemma} 
\begin{remark} By using Gagliardo-Nirenberg inequality in Lemma
\ref{lemma2.1}, let $p=3,q=r=2,\lambda^{*}=\frac{N}{6}$, there exists constant
$C_{GN}>0$, such that \begin{equation}\label{22}
  \int_{\mathbb{R}^{N}}u^{3}dy\leq C_{GN}(N)(\Vert\nabla u\Vert_{L^{2}(\mathbb{R}^{N})}^{\frac{N}{2}}\Vert u\Vert_{L^{2}(\mathbb{R}^{N})}^{3-\frac{N}{2}}+\Vert u\Vert_{L^{2}(\mathbb{R}^{N})}^{3}).
\end{equation} 
On the one hand, by Young's inequality, we can obtain
\begin{align*} 
&C_{GN}(N)\Vert\nabla
u\Vert_{L^{2}(\mathbb{R}^{N})}^{\frac{N}{2}}\Vert
u\Vert_{L^{2}(\mathbb{R}^{N})}^{3-\frac{N}{2}}\\
&\leq \frac{1}{\mu}\Vert \nabla
u\Vert_{L^{2}(\mathbb{R}^{N})}^{2}+\mu^{\frac{N}{4-N}}C_{GN}^{\frac{4}{4-N}}(N)\Vert
u\Vert_{L^{2}(\mathbb{R}^{N})}^{\frac{2(6-N)}{4-N}}. \end{align*} and
\begin{equation}\label{24}
  C_{GN}(N)\left \| u \right \|_{L^{2}(\mathbb{R}^{N})}^{3}\leq \left \| u \right \|_{L^{2}(\mathbb{R}^{N})}^{\frac{2(6-N)}{4-N}}+C_{GN}^{\frac{2(6-N)}{N}}(N).
\end{equation} By interpolation inequality, we obtain
\begin{equation}\label{25} \begin{aligned} \left \| u \right
\|_{L^{2}(\mathbb{R}^{N})}^{\frac{2(6-N)}{4-N}}&\leq (\left \| u \right
\|_{L^{1}(\mathbb{R}^{N})}^{\frac{1}{4}}\left \| u \right
\|_{L^{3}(\mathbb{R}^{N})}^{\frac{3}{4}})^{\frac{2(6-N)}{4-N}}\\ &=(\int
_{\mathbb{R}^{N}}u^{3}dy\int _{\mathbb{R}^{N}}udy)^{\frac{6-N}{2(4-N)}}.
\end{aligned} \end{equation} \end{remark}

\begin{lemma}\textsuperscript{\cite{bourgain2001limiting}}
  The embedding $W^{s,p}(\mathbb{R}^{N})\hookrightarrow L^{r}(\mathbb{R}^{N})$ is continuous, that is \begin{equation}\label{55yy}
      \left \| u \right \|_{L^{r}}\leq C_{*}\left [ u \right ]_{W^{s,p}},\quad \forall u\in W_{0}
  \end{equation}
  where $1\leq r<p_{s}^{*},C_{*}=C_{*}(s,p,r,N,\mathbb{R}^{N})$ is optimal embedded constant. And $W_{0}$ is defined
  $$W_{0}:=\left \{ u\in W^{s,p}(\mathbb{R}^{N}) |  \int_{\mathbb{R}^{N}}udx=0\right \}.$$
\end{lemma}

\begin{lemma}\textsuperscript{\cite{2015Maximum}}\label{7700}
   Suppose that a nonnegative function $y_{k}(t)\geq 0, \tilde{a},\beta >0$  be a solution of the fractional differential inequality
   \begin{equation}\label{2655}
       _{0}^{C}\textrm{D}_{t }^{\alpha }y_{k}(t)\leq -\tilde{a}y_{k}(t)+\beta^{-1}b_{k}(t)y_{k}^{1-\beta}(t).
   \end{equation}
Then, the solution of \eqref{2655} can be estimated as \begin{align*}
    y_{k}(t)&\leq y_{k}(0)+\lambda _{k}y_{k}(0)\int_{0}^{t}(t-s)^{\alpha-1}E_{\alpha ,\alpha }(\lambda _{k}(t-s)^{\alpha })ds\\
    &\qquad+\varepsilon ^{\frac{1}{\beta }}\int_{0}^{t}(t-s)^{\alpha-1}E_{\alpha ,\alpha }(\lambda _{k}(t-s)^{\alpha })b_{k}^{\frac{1}{\beta }}(s)ds,
\end{align*} where $\lambda _{k}=-\tilde{a}+\frac{1-\beta }{\beta \varepsilon
^{\frac{1}{1-\beta }}}$ and $\varepsilon >0$ \end{lemma}

\begin{definition}\textsuperscript{\cite{2006Global}}
  Let $X$ be a Banach space,$z_{0}$ belong to $X$, and $f\in L^{1}(0,T;X)$. The function $z(x,t)\in C([0,T];X)$ given by
  \begin{equation}\label{pp}
      z(t)=e^{-t}e^{t\Delta}z_{0}+\int_{0}^{t}e^{-(t-s)}\cdot e^{(t-s)\Delta }f(s)ds, 0\leq t\leq T,
  \end{equation}
  is the mild solution of \eqref{pp} on $[0,T]$, where $(e^{t\Delta}f)(x,t)=\int _{\mathbb{R}^{N}}G(x-y,t)f(y)dy$ and $G(x,t)$ is the heat kernel by $G(x,t)=\frac{1}{(4\pi t)^{N/2}}exp(-\frac{\left | x \right |^{2}}{4t})$.
\end{definition} 
\begin{lemma}\textsuperscript{\cite{2006Global}}\label{lemma5}
  Let $0\leq q\leq p\leq \infty ,\frac{1}{q}-\frac{1}{p}<\frac{1}{N}$ and suppose that $z$ is the function given by \eqref{pp} and $z_{0}\in W^{1,p}(\mathbb{R}^{N})$. If $f\in L^{\infty }(0,\infty ;L^{q}(\mathbb{R}^{N}))$, then
  \begin{eqnarray}
   \nonumber \left \| z(t) \right \|_{L^{p}(\mathbb{R}^{N})}&\leq& \left \|  z_{0}\right \|_{L^{p}(\mathbb{R}^{N})}+C\cdot \Gamma (\gamma )\underset{0<s<t}{sup}\left \| f(s) \right \|_{L^{q}(\mathbb{R}^{N})},\\
\nonumber\left \| \nabla z(t) \right \|_{L^{p}(\mathbb{R}^{N})}&\leq& \left \|
\nabla z_{0}\right \|_{L^{p}(\mathbb{R}^{N})}+C\cdot \Gamma (\tilde{\gamma
})\underset{0<s<t}{sup}\left \| f(s) \right \|_{L^{q}(\mathbb{R}^{N})},
  \end{eqnarray}
  for $t\in [0,\infty )$, where $C$ is a positive constant independent of $p,\Gamma (\cdot )$ is the gamma function, and $\gamma =1-(\frac{1}{q}-\frac{1}{p})\cdot \frac{N}{2},\tilde{\gamma }=\frac{1}{2}-(\frac{1}{q}-\frac{1}{p})\cdot \frac{N}{2}$.
\end{lemma}

\begin{lemma}\textsuperscript{\cite{2016fractional}}\label{233}
  Assume that $(a,b)\in (\mathbb{R}^{+})^{2},0<\alpha <1$, then there exist $c_{1},c_{2},c_{3}>0$, such that
  \begin{equation*}
      (a+b)^{\alpha }\leq c_{1}a^{\alpha }+c_{2}b^{\alpha }
  \end{equation*}
  and
  \begin{equation}\label{2.3}
      \left | a-b \right |^{p-2}(a-b)(a^{\alpha }-b^{\alpha })\geq c_{3}\left | a^{\frac{p+\alpha-1 }{p}}-b^{\frac{p+\alpha-1 }{p}} \right |^{p}.
  \end{equation}
\end{lemma}

\begin{lemma}\textsuperscript{\cite{2014Ultra}}\label{5566}
  Let $N\geq 3,q>1,m>1-2/N$, assume $u\in L_{+}^{1}(\mathbb{R}^{N})$ and $u^{\frac{m+q-1}{2}}\in H^{1}(\mathbb{R}^{N})$, then
  \begin{equation*}
      (\left \| u \right \|_{q}^{q})^{1+\frac{m-1+2/N}{q-1}}\leq S_{N}^{-1}\left \| \nabla u^{(q+m-1)/2} \right \|_{2}^{2}\left \| u \right \|_{1}^{\frac{1}{q-1}(2q/N+m-1)}.
  \end{equation*}
\end{lemma}

 \begin{lemma}\textsuperscript{\cite{1959On}}\textsuperscript{\cite{1971Propriet}}\label{lemma33}
   When the parameters $p,q,r$ meet any of the following conditions:

   (i) $q >N \geq 1$, $r\geq1$ and $p=\infty$;

   (ii)$q>max\left \{ 1,\frac{2N}{N+2} \right \},1\leq r<\sigma $ and $r<p<\sigma +1$ in
   \begin{equation}
    \nonumber \sigma :=\begin{cases}
               \frac{(q-1)N+q}{N-q},&q<N,\\
               \infty ,&q\geq N.
              \end{cases}
   \end{equation}
Then the following inequality is established $$\left \| u \right
\|_{L^{p}(\mathbb{R}^{N})}\leq C_{GN}\left \| u \right
\|_{L^{r}(\mathbb{R}^{N})}^{1-\lambda ^{*}} \left \| \nabla u \right
\|_{L^{q}(\mathbb{R}^{N})}^{\lambda ^{*}}$$ among $$\lambda
^{*}=\frac{qN(p-r)}{p[N(q-r)+qr]}.$$
 \end{lemma}
 \begin{remark}
     The following estimate $k\int _{\mathbb{R}^{N}}u^{k+1}dx$. when $m\leq3$ and
\begin{equation}\label{5.5}
 k> max\left \{ \frac{(3-m)(N-2)}{4} -(m-1),(\frac{2-m}{2})N-1,N(2-m)-2\right \},
\end{equation} from Lemma \ref{lemma33}, we obtain $$ \begin{aligned} &k\int
_{\mathbb{R}^{N}}u^{k+1}dx=k\left \| u^{\frac{k+m-1}{2}} \right
\|_{L^{\frac{2(k+1)}{k+m+1}}(\mathbb{R}^{N})}^{\frac{2(k+1)}{k+m-1}}\\ &\leq
k\left \| u^{\frac{k+m-1}{2}} \right
\|_{L^{\frac{2(k+1)}{k+m+1}}(\mathbb{R}^{N})}^{\frac{2(k+1)}{k+m-1}\lambda
^{*}}\left \| \nabla u^{\frac{k+m-1}{2}} \right
\|_{L^{2}(\mathbb{R}^{N})}^{\frac{2(k+1)(1-\lambda ^{*})}{k+m-1}}.
\end{aligned} $$ where $$\lambda
^{*}=\frac{1+\frac{(k+m-1)N}{2(k+1)}-\frac{N}{2}}{1+\frac{(k+m-1)N}{k+2}-\frac{N}{2}}\in
\left \{ max\left \{ 0,\frac{2-m}{k+1} \right \},1 \right \},$$ using Young's
inequality, there are \begin{equation}\label{4.1.3} \begin{aligned} k\int
_{\mathbb{R}^{N}}u^{k+1}dx&\leq k\left \| u^{\frac{k+m-1}{2}} \right
\|_{L^{\frac{k+2}{k+m+1}}(\mathbb{R}^{N})}^{\frac{2(k+1)}{k+m-1}\lambda
^{*}}\left \| \nabla u^{\frac{k+m-1}{2}} \right
\|_{L^{2}(R^{N})}^{\frac{2(k+1)(1-\lambda ^{*})}{k+m-1}}\\ &\leq
\frac{2k(k-1)}{(m+k-1)^{2}}\left \| \nabla u^{\frac{k+m-1}{2}} \right
\|_{L^{2}(R^{N})}^{2}\\ &\qquad+C_{1}(N,k,m)\left \| u^{\frac{k+m-1}{2}} \right
\|_{L^{\frac{k+2}{k+m-1}}(\mathbb{R}^{N})}^{Q_{2}} \end{aligned} \end{equation}
where is $Q_{2}=\frac{2(k+1)\lambda ^{*}}{m-2+(k+1)\lambda ^{*}}$. Next
estimate $\left \| u^{\frac{k+m-1}{2}} \right
\|_{L^{\frac{k+2}{k+m-1}}(\mathbb{R}^{N})}^{Q_{2}}$. We will use the
interpolation inequality to get \begin{equation}\label{4.1.4}
  \begin{aligned}
  &\left \| u^{\frac{k+m-1}{2}} \right \|_{L^{\frac{k+2}{k+m-1}}(\mathbb{R}^{N})}^{Q_{2}}\leq\left \| u^{\frac{k+m-1}{2}} \right \|_{L^{\frac{2(k+1)}{k+m-1}}(\mathbb{R}^{N})}^{Q_{2}\lambda }\left \| u^{\frac{k+m-1}{2}} \right \|_{L^{\frac{2}{k+m-1}}(\mathbb{R}^{N})}^{Q_{2}(1-\lambda )}\\
&\leq(\left \| u^{\frac{k+m-1}{2}} \right
\|_{L^{\frac{2(k+1)}{k+m-1}}(\mathbb{R}^{N})}^{\frac{2(k+1)}{k+m-1} }\left \|
u^{\frac{k+m-1}{2}} \right
\|_{L^{\frac{2}{k+m-1}}(\mathbb{R}^{N})}^{\frac{2}{k+m-1}})^{\frac{Q_{2}\lambda
(k+m-1)}{2(k+1)}}\\ &\left \| u^{\frac{k+m-1}{2}} \right
\|_{L^{\frac{2}{k+m-1}}(\mathbb{R}^{N})}^{Q_{2}(1-\lambda -\frac{\lambda
}{k+1})}
  \end{aligned}
\end{equation} in $\lambda =\frac{k+1}{k+2}$, and $$Q_{2}(1-\lambda
-\frac{\lambda }{k+1})=0.$$ Then 
\begin{equation}\label{4.1.5} \begin{aligned}
  &C_{1}(N,k,m)\left \| u^{\frac{k+m-1}{2}} \right \|_{L^{\frac{k+2}{k+m-1}}(\mathbb{R}^{N})}^{Q_{2}}\\
  &\leq C_{1}(N,k,m)\\
  &\,(\left \| u^{\frac{k+m-1}{2}} \right \|_{L^{\frac{2(k+1)}{k+m-1}}(\mathbb{R}^{N})}^{\frac{2(k+1)}{k+m-1} }\left \| u^{\frac{k+m-1}{2}} \right \|_{L^{\frac{2}{k+m-1}}(R^{N})}^{\frac{2}{k+m-1}})^{\frac{Q_{2}\lambda (k+m-1)}{2(k+1)}}.
\end{aligned} \end{equation} Noticing that when $m>2-\frac{2}{N}$, it is easy
to verify$$\frac{Q_{2}\lambda (k+m-1)}{2(k+1)}=\frac{(k+1)(m-2)+\lambda
^{*}(k+1)(3-m)}{(k+1)(k+m-1)\lambda ^{*}}< 1,$$ using Young's inequality, then
\begin{equation}\label{4.1.6} \begin{aligned} &C_{1}(N,k,m)(\left \|
u^{\frac{k+m-1}{2}} \right
\|_{L^{\frac{2(k+1)}{k+m-1}}(\mathbb{R}^{N})}^{\frac{2(k+1)}{k+m-1} }\left \|
u^{\frac{k+m-1}{2}} \right
\|_{L^{\frac{2}{k+m-1}}(\mathbb{R}^{N})}^{\frac{2}{k+m-1}})^{\frac{Q_{2}\lambda
(k+m-1)}{2(k+1)}}\\
 &\leq k\left \| u^{\frac{k+m-1}{2}} \right \|_{L^{\frac{2(k+1)}{k+m-1}}(\mathbb{R}^{N})}^{\frac{2(k+1)}{k+m-1} }\left \| u^{\frac{k+m-1}{2}} \right \|_{L^{\frac{2}{k+m-1}}(\mathbb{R}^{N})}^{\frac{2}{k+m-1}}\\
 &\qquad+C_{2}(N,k,m).
\end{aligned} \end{equation}
 \end{remark}
\begin{lemma}\textsuperscript{\cite{Shen2016A}}\label{lemma333}
  Let $N\geq 1$. $p$ is the exponent from the Sobolev embedding theorem, $i.e.$
  \begin{equation}\label{22.1}
    \begin{cases}
      p=\frac{2N}{N-2},& N\geq 3\\
     2<p<\infty,& N=2\\
    p=\infty,& N=1
    \end{cases}
  \end{equation}
  $1\leq r<q<p$ and $\frac{q}{r}< \frac{2}{r}+1-\frac{2}{p}$, then for $v\in {H}'(\mathbb{R}^{N})$ and $v\in L^{r}(\mathbb{R}^{N})$, it holds
  \begin{equation}\label{99.1}
  \begin{aligned}
    \left \| v \right \|_{L^{q}(\mathbb{R}^{N})}^{q}&\leq
     C(N)c^{\frac{\lambda q}{2-\lambda q}}_{0}\left \| v \right \|_{L^{p}(\mathbb{R}^{N})}^{\gamma }+c_{0}\left \| \nabla u \right \|_{L^{2}(\mathbb{R}^{N})}^{2},\quad N>2,\\
\left \| v \right \|_{L^{q}(\mathbb{R}^{N})}^{q}&\leq C(N)(c^{\frac{\lambda
q}{2-\lambda q}}_{0}+c^{-\frac{\lambda q}{2-\lambda q}}_{1})\left \| v \right
\|_{L^{r}(\mathbb{R}^{N})}^{\gamma }+c_{0}\left \| \nabla u \right
\|_{L^{2}(\mathbb{R}^{N})}^{2}\\ &\qquad+c_{1}\left \| v \right
\|_{L^{2}(\mathbb{R}^{N})}^{2},\quad  N=1,2.
	\end{aligned}
  \end{equation}
  Here $C(N)$ are constants depending on $N,c_{0},c_{1}$ are arbitrary positive constants and
  \begin{equation}\label{112}
  \nonumber \lambda =\frac{\frac{1}{r}-\frac{1}{q}}{\frac{1}{r}-\frac{1}{p}}\in (0,1),\gamma =\frac{2(1-\lambda )q}{2-\lambda q}=\frac{2(1-\frac{q}{p})}{\frac{2-q}{r}-\frac{2}{p}+1}.
  \end{equation}
\end{lemma} 
\begin{remark} By Lemma \ref{lemma333}, letting $$v=\frac{m+q_
{k}-1}{2},q=\frac{2(q_{k}+1)}{m+q_{k}-1},r=\frac{2q_{k-1}}{m+q_{k}-1},c_{0}=c_{1}=\frac{1}{2q_{k}},$$
one has that for $N \geq 3$, 
\begin{equation}\label{4.2.2} 
\begin{aligned}
  &\left \| u \right \|_{L^{q_{k}+1}(\mathbb{R}^{N})}^{q_{k}+1}\leq C(N)c_{0}^{\frac{1}{\delta _{1}-1}}(\int _{\mathbb{R}^{N}}u^{q_{k-1}}dx)^{\gamma _{1}}
\\ &\qquad+\frac{1}{2q^{k}}\left \| \nabla  u^{\frac{m+q_{k}-1}{2}} \right
\|_{L^{2}(\mathbb{R}^{N})}^{2}+\frac{1}{2q_{k}}\left \| u \right
\|_{L^{m+q_{k}-1}(\mathbb{R}^{N})}^{m+q_{k}-1},
\end{aligned}
\end{equation}
where $$\gamma _{1}=1+\frac{q_{k}+q_{k-1}+1}{q_{k-1}+\frac{p(m-2)}{p-2}}\leq
2,$$ $$\delta
_{1}=\frac{(m+q_{k}-1)-2\frac{q_{k-1}}{p^{*}}}{q_{k}-q_{k-1}+1}=O(1).$$

Applying Lemma \ref{lemma333} with $$v=u^
{\frac{m+q_{k}-1}{2}},q=2,r=\frac{2q_{k-1}}{m+q_{k}-1},c_{0}=c_{1}=\frac{1}{2}$$
noticing $q_{k-1}=\frac{(q_{k}+1)+1}{2}$, and using Young's inequality, we
obtain \begin{equation}\label{4.2.4}
  \begin{aligned}
  &\frac{1}{2}\left \| u \right \|_{m+q_{k}-1}^{m+q_{k}-1}=\frac{1}{2}\int _{\mathbb{R}^{N}}u^{m+q_{k}-1}dx\\
&\leq c_{2}(N)(\int _{\mathbb{R}^{N}}u^{q_{k-1}}dx)^{\gamma
_{2}}+\frac{1}{2}\left \| \nabla u^{\frac{m+q_{k}-1}{2}} \right
\|_{L^{2}(\mathbb{R}^{N})}^{2}\\ &\leq \int _{\mathbb{R}^{N}}udx\int
_{\mathbb{R}^{N}}u^{q_{k}+1}dx+c_{3}(N)\\
&\qquad+\frac{1}{2}\left \| \nabla
u^{\frac{m+q_{k}-1}{2}} \right \|_{L^{2}(\mathbb{R}^{N})}^{2},
  \end{aligned}
\end{equation} where $$\gamma _{2}=1+\frac{m+q_{k}-q_{k-1}-1}{q_{k-1}}< 2.$$
\end{remark}
\begin{lemma}\label{ollkkkmh} 
(Comparison principle)
Let $0<s<1,p>1,\Omega \subset \mathbb{R}^{N}$ and let $u,v\in \Pi$ be a real-valued weak subsolution and
supersolution of \eqref{1}-\eqref{2}, respectively, with $u_{0}(x)\leq
v_{0}(x)$ for $x\in \Omega$. Based on the following fractional differential inequality
\begin{align*}
 \nonumber   &\int_{0}^{t}\int _{\Omega }(_{0}^{C}\textrm{D}_{\tau }^{\alpha}[u-v])(u-v)_{+} dxd\tau\\
   & \leq \left ( \mu L(2)-k\mu \eta ^{\frac{1}{2}}C_{4}-\gamma  \right )\int_{0}^{t}\int _{\Omega }(u-v)(u-v)_{+} dxd\tau,
\end{align*} 
then $u\leq v$ a.e. in $\Omega_{T}$, where
$\Omega_{T}=\Omega\times(0,T)$ and
$L(m)=C(m)max(\left \| u \right \|_{L^{\infty
}(\Omega)}^{m-1},\left \| v \right \|_{L^{\infty }(\Omega)}^{m-1}).$
\end{lemma} 
\begin{pro}
 We select the function $\varphi=(u-v)_{+}$, where $\varphi=(u-v)_{+}$ is a positive real part of the real number $\varphi=(u-v)_{+}=max\left \{u-v,0 \right \}$. Followed by  $\varphi(x,0)=0,\varphi (x,t)\mid_{\partial \Omega }=0$, then we obtain for $t\in (0,T]$
 \begin{align}\label{fiofio}
\nonumber     &\int_{0}^{t}\int _{\Omega }\,_{0}^{C}\textrm{D}_{\tau
}^{\alpha}[u-v]\varphi dxd\tau+\int_{0}^{t}\int _{\Omega }[(-\Delta
)_{p}^{s}u-(-\Delta )_{p}^{s}v]\varphi dxd\tau \\ \nonumber&\leq \mu
\int_{0}^{t}\int _{\Omega }(u^{2}-v^{2})\varphi dxd\tau -k\mu\int_{0}^{t}\int
_{\Omega }(u^{2}J*u-v^{2}J*v)\varphi dxd\tau\\ &\qquad-\gamma \int_{0}^{t}\int
_{\Omega }(u-v)\varphi dxd\tau
 \end{align}
 By \eqref {4}, we can write the last term on the left of the inequality \eqref {fiofio} as
\begin{align}\label{onhu}
 \nonumber   &\int_{0}^{t}\int _{\Omega }[(-\Delta )_{p}^{s}u-(-\Delta )_{p}^{s}v]\varphi dxd\tau\\
\nonumber&=\int_{0}^{t}\int _{\Omega }\int _{\Omega }\frac{\left | u(x)-u(y)
\right |^{p-2}(u(x)-u(y))(\varphi (x)-\varphi (y))}{\left | x-y \right
|^{N+sp}}\\
\nonumber&\qquad dxdyd\tau \\ 
\nonumber&-\int_{0}^{t}\int _{\Omega }\int _{\Omega
}\frac{\left | v(x)-v(y) \right |^{p-2}(v(x)-v(y))(\varphi (x)-\varphi
(y))}{\left | x-y \right |^{N+sp}}\\
\nonumber&\qquad dxdyd\tau \\ &=\int_{0}^{t}\int _{\Omega
}\int _{\Omega }\frac{\mathcal{M}(u,v)(\varphi (x)-\varphi (y))}{\left | x-y
\right |^{N+sp}}dxdyd\tau , 
\end{align} where
\begin{align}\label{5ghgh}
\nonumber    &\mathcal{M}(u,v)=\left | u(x)-u(y) \right |^{p-2}(u(x)-u(y))\\
    &\quad-\left | v(x)-v(y) \right |^{p-2}(v(x)-v(y))(\varphi (x)-\varphi (y))
\end{align} 
Thus, we can show that \begin{align*}
    &\mathcal{M}(u,v)(\varphi (x)-\varphi (y))=[\left | u(x)-u(y) \right |^{p-2}(u(x)-u(y))\\
    &-\left | v(x)-v(y) \right |^{p-2}(v(x)-v(y))(\varphi (x)-\varphi (y))]\\
    &\times [(u(x)-u(y))-(v(x)-v(y))]_{+}
\end{align*} is nonnegative for any $p>1$, thanks to the two inequalitys
(see\cite{lindqvist2019notes},P,99-100) 
\begin{equation*} 
 \left |
\left | a \right |^{\frac{p-2}{2}}a-\left | b \right |^{\frac{p-2}{2}}b  \right
|^{2}\leq \left \langle \left | a \right |^{p-2}a-\left | b \right |^{p-2}b,a-b
\right \rangle,\quad p\geq 2
\end{equation*}
and
\begin{align*}
   & (p-1)\left | b-a \right |^{2}\int_{0}^{1}\left
| a+t(b-a) \right |^{p-2}dt\\
&\leq \left \langle \left | a \right |^{p-2}a-\left |
b \right |^{p-2}b,a-b \right \rangle,\,  1\leq p\leq 2 
\end{align*}
with
$a:=u(x)-u(y),b:=v(x)-v(y)$ in \eqref{5ghgh}.
Now, we'll zoom in on the right side of \eqref {fiofio}.

By [ \cite{borikhanov2022qualitative},Theorem 3.2
p12], we take into account the following inequalities for $m\geq 2$ 
\begin{align*}
    \left | \left | u
\right |^{m-1}u-\left | v \right |^{m-1}v \right |&\leq C(m)\left | u-v \right
|\left | \left | u \right |^{m-1}+\left | v \right |^{m-1} \right |\\
&\leq
L(m)\left | u-v \right |,
\end{align*}
where $L(m)=C(m)max(\left \| u \right \|_{L^{\infty
}(\Omega)}^{m-1},\left \| v \right \|_{L^{\infty }(\Omega)}^{m-1})$. Specially, when $m=2$, then \begin{equation}\label{380}
   \int_{0}^{t}\int _{\Omega }(u^{2}-v^{2})\varphi dxd\tau \leq L(2)\int_{0}^{t}\int _{\Omega } (u-v)\varphi dxd\tau,
\end{equation} 
where we used the commonly used inequality 
[\cite{2012Hitchhiker}, Theorem 8.2] for any $u\in L^{p}(\Omega)$ suth that
$$\left \| u \right \|_{C(\Omega )}\leq \left \| u \right \|_{C^{0,\beta
}(\Omega )}\leq \left \| u \right \|_{W^{s,p}(\Omega )},\beta =(sp-N)/p ,$$
which gives the boundness of $max(\left \| u \right \|_{C(\Omega
)(\Omega)}^{m-1},\left \| v \right \|_{C(\Omega )(\Omega)}^{m-1})$.

In addition, by Theorem \ref{theorem4}, let $$C_{4}=min \left \{ \left \| u
\right \|_{L^{1}(\Omega )}^{\frac{1}{2}},\left \| v \right \|_{L^{1}(\Omega
)}^{\frac{1}{2}} \right \},$$ then we have \begin{align}\label{3119}
 \nonumber  & k\mu\int_{0}^{t}\int _{\Omega }(u^{2}J*u-v^{2}J*v)\varphi dxd\tau\\
\nonumber&=k\mu\int_{0}^{t}\int _{\Omega }(u(\int_{\Omega
}J(x-y)u(y,t)dy)^{\frac{1}{2}}-v(\int_{\Omega }J(x-y)v(y,t)dy)^{\frac{1}{2}})\\
\nonumber&\times (u(\int_{\Omega }J(x-y)u(y,t)dy)^{\frac{1}{2}}+v(\int_{\Omega
}J(x-y)v(y,t)dy)^{\frac{1}{2}})dxd\tau \\ \nonumber&\geq k\mu\int_{0}^{t}\int
_{\Omega }(u(\int_{\Omega }J(x-y)u(y,t)dy)^{\frac{1}{2}}-v(\int_{\Omega
}J(x-y)v(y,t)dy)^{\frac{1}{2}})\\
\nonumber&\qquad dxd\tau\\
\nonumber&\geq k\mu\eta
^{\frac{1}{2}}\int_{0}^{t}\int _{\Omega }(u(\int_{\Omega
}u(y,t)dy)^{\frac{1}{2}}-v(\int_{\Omega }v(y,t)dy)^{\frac{1}{2}})dxd\tau\\
&\geq k\mu\eta ^{\frac{1}{2}}C_{4}\int_{0}^{t}\int _{\Omega }(u-v)\varphi
dxd\tau. \end{align} Combining \eqref{onhu},\eqref{380} and \eqref{3119}, we
can rewrite the inequality \eqref{fiofio} as \begin{align}\label{gkdf}
 \nonumber   &\int_{0}^{t}\int _{\Omega }(_{0}^{C}\textrm{D}_{\tau }^{\alpha}[u-v])(u-v)_{+} dxd\tau\\
   & \leq \left ( \mu L(2)-k\mu \eta ^{\frac{1}{2}}C_{4}-\gamma  \right )\int_{0}^{t}\int _{\Omega }(u-v)(u-v)_{+} dxd\tau.
\end{align} 
Using the lemma\ref {lemma1}, we can rewrite the inequality \eqref {gkdf} as follows 
\begin{align}\label{31489}
 \nonumber   &\frac{1}{2}\int_{0}^{t}\int _{\Omega }\,_{0}^{C}\textrm{D}_{\tau }^{\alpha}(u-v)_{+}^{2} dxd\tau\\
   & \leq \left ( \mu L(2)-k\mu \eta ^{\frac{1}{2}}C_{4}-\gamma  \right )\int_{0}^{t}\int _{\Omega }(u-v)_{+}^{2} dxd\tau.
\end{align} 
Because $k,\mu,\eta>0,\gamma \geq 1$, applying left Caputo
fractional differential operator $_{0}^{C}\textrm{D}_{t }^{1-\alpha}$  to both sides of\eqref {31489} and using the lemma\ref {gugui}, we get
\begin{equation*}
    \frac{1}{2}\int _{\Omega }(u-v)_{+}^{2} dx\leq \mu L(2)\int _{\Omega }\int_{0}^{t}(t-\tau)^{\alpha-1}(u-v)_{+}^{2} d\tau dx.
\end{equation*}
Then, from the weakly singular Gronwall's inequality [see
\cite{ball1982geometric} Lemma 7.1.1] \begin{equation*}
    \int _{\Omega }(u-v)_{+}^{2}dx=0\Leftrightarrow (u-v)_{+}=0,\quad x\in \Omega.
\end{equation*} Finally, using the fact $(u-v)_{+}=max\left \{ u-v,0 \right
\}$, it follows that $u\leq v$ almost everywhere for $(x,t)\in \Omega_{T}$.
\end{pro} \begin{remark} The method of proving this Lemma is similar to
[\cite{borikhanov2022qualitative}, Theorem 3.2]. But the biggest difference
between the two paper is that the range of $(-\Delta)_{p}^{s},p$ value is
different. This paper is $(-\Delta)_{p}^{s},p>1$, and the literature
\cite{borikhanov2022qualitative} is $p\geq2$. Not only that, but also slightly
different about the a shrinking part of the non-local items for space-time
fractional diffusion equation. \end{remark}

\section{Existence and unique
weak solution for \eqref{114}-\eqref{1.1.5} } 
\subsection{Weighted $L^{1}$ estimates in the nonlinear fractional diffusion
range}\label{zhanhui22}
  In order to state the main estimate we introduce that concept of weighted mass at time
$0\leq t<T$ \begin{equation}\label{mkmk}
    X(t;u,v,\varphi )=\int _{\mathbb{R}^{N}}(u(t)-v(t))\varphi dx,
\end{equation} where the weight $\varphi$ is a positive function to be
specified next. From \cite{vazquez2022growing}, we introduce the operator
$\mathcal{M}_{{s}'}$ by the formula \begin{equation}
    (\mathcal{M}_{{s}'}\varphi)(x):=P.V.\int _{\mathbb{R}^{N}}\frac{\left | \varphi (x,t)-\varphi (y,t)  \right |}{\left | x-y \right |^{N+2{s}'}}dy.
\end{equation} We remark that when $0<2{s}'<1$ this operator is well-defined
and bounded for bounded and uniformly Lipschitz continuous function since the
singularity at $x=y$ is integrable.

\textbf{The class $\mathcal{C}=\mathcal{C}(s,p,m)$.}\quad The class of suitable
weight function for our main estimate is formed by the smooth and positive
functions $\varphi$ defined in $\mathbb{R}^{N}$ such that
$\mathcal{M}_{sp/2}(\varphi)$ is locally bounded and
\begin{equation}\label{opop}
    C(\varphi )=\int _{\mathbb{R}^{N}}\frac{\left |\mathcal{M}_{sp/2}\varphi(x) \right |^{\frac{1}{1-m(p-1)}}}{\varphi (x)^{\frac{m(p-1)}{1-m(p-1)}}}dx< \infty.
\end{equation} Note that this class depends on $s,p,m$. The conditions
${s}'=sp/2<1$ ensures that the class contains a large class of uniformly
Lipschitz functions depending on our choice of $s,p$ and $m$. The value of
$C(\varphi )$ only depends on the positivity, smoothness and behaviour of
$\varphi (x)$ as $\left | x \right |\rightarrow \infty $.

\textbf{Admissible decay rates.}\quad From \cite{bonforte2014quantitative},
there are many smooth, bounded and positive function $\varphi$ decaying at
infinity like a power $\varphi \sim O(\left | x \right |^{-(N+\gamma )})$ with
$\gamma >2{s}'$ such that $(-\Delta )^{{s}'}\varphi $ decays like $O(\left | x
\right |^{-(N+2{s}')})$ as $\left | x \right |\rightarrow \infty $, assuming
that $0<{s}'<1$. We can check that $\mathcal{M}_{{s}'}\varphi $ decays in the
same way if $2{s}'<1$. We have $$\frac{\left |\mathcal{M}_{sp/2}\varphi(x)
\right |^{\frac{1}{1-m(p-1)}}}{\varphi (x)^{\frac{m(p-1)}{1-m(p-1)}}}dx\sim
\left | x \right |^{-\mu }$$ with $$\mu
=\frac{(N+2{s}')}{1-m(p-1)}-\frac{(N+\gamma )(m(p-1))}{1-m(p-1)}.$$ The
expression is integrable if $\mu>N$. Working out the details we find that
$C(\varphi)$ is finite if $\gamma <\frac{2{s}'}{m(p-1)}.$
\begin{thm} (Weighted $L^{1}$ estimates).\quad Let
$0<s<1,1<p<2,0<m<\frac{1}{p-1},\varepsilon>0$ with $sp<1$. Let $u\geq v$ be two
nonnegative semigroup solution of \eqref{100} in a strip
$Q_{T}=\mathbb{R}^{N}\times (0,T)$ with $T>0$ and $X(t)$ be as in \eqref{mkmk}.
Then for all $\varphi \in \mathcal{C}(s,p,m)$ there is a finite constant $K>0$
depending on $\varphi,\varepsilon $ such that we have
\begin{align}\label{dfdf}
 \nonumber  & \left | X^{1-m(p-1)}(t_{1})- X^{1-m(p-1)}(t_{2}) \right |\\
    &\leq K(\varphi)\left | t_{1}-t_{2} \right |^{\alpha [1-m(p-1)] }.
\end{align} Actually, we may take $K(\varphi)=\varepsilon \left( \frac{
C(\varphi )}{\alpha \Gamma (\alpha )}\right )^{1-m(p-1)}$ with $C(\varphi )$
given by \eqref{opop}. \end{thm} \begin{pro} We multiply an $L^{2}$ solution by
a smooth and positive test function $\varphi (x),u\geq v$, and from Lemma
\ref{lemma1}, we have 
\begin{align*}
   & \left |\int _{\mathbb{R}^{N}} (_{0}^{C}\textrm{D}_{t}^{\alpha }u(x,t)-_{0}^{C}\textrm{D}_{t}^{\alpha }v(x,t))\varphi (x)dx \right |\\
&=\left | _{0}^{C}\textrm{D}_{t}^{\alpha }\int
_{\mathbb{R}^{N}}(u(x,t)-v(x,t))\varphi (x)dx \right |\\ &\leqslant\left | \int
_{\mathbb{R}^{N}}((-\Delta )_{p}^{s}u^{m}-(-\Delta )_{p}^{s}v^{m})\varphi (x)dx
\right |+\left | \int _{\mathbb{R}^{N}}(u^{2}-v^{2})\varphi dx \right |\\
&\qquad-\int _{\mathbb{R}^{N}}\left [ u^{2}\int _{\mathbb{R}^{N}}udy-v^{2}\int
_{\mathbb{R}^{N}}vdy \right ]\varphi dx-\int _{\mathbb{R}^{N}}(u-v)\varphi
dx.\\ &\leqslant\left | \int _{\mathbb{R}^{N}}(-\Delta
)_{p}^{s}(u^{m}(x,t)-v^{m}(x,t))\varphi (x)dx \right |+\left | \int
_{\mathbb{R}^{N}}(u^{2}-v^{2})\varphi dx \right |\\
&=\iint_{\mathbb{R}^{2N}}\frac{\left [
(u^{m}(x,t)-v^{m}(x,t))^{p-1}-(u^{m}(y,t)-v^{m}(y,t))^{p-1} \right ]}{\left |
x-y \right |^{N+sp}}\\ &\qquad(\varphi (x)-\varphi (y))dxdy+\left | \int
_{\mathbb{R}^{N}}(u^{2}-v^{2})\varphi dx \right |\\ 
&\leqslant 2\int
_{\mathbb{R}^{N}}(u^{m}(x,t)-v^{m}(x,t))^{p-1}dx\left ( \int
_{\mathbb{R}^{N}}\frac{\left | \varphi (x)-\varphi (y) \right |}{\left | x-y
\right |^{N+sp}}dy \right )\\
&\quad+\left | \int _{\mathbb{R}^{N}}(u^{2}-v^{2})\varphi
dx \right |\\
&=2\int
_{\mathbb{R}^{N}}(u^{m}(x,t)-v^{m}(x,t))^{p-1}\mathcal{M}_{{s}'}\varphi(x)
dx+\left | \int _{\mathbb{R}^{N}}(u^{2}-v^{2})\varphi dx \right |\\
&\leqslant_{(a)}2\int _{\mathbb{R}^{N}}\left [ (u-v)^{m} \right
]^{p-1}\mathcal{M}_{{s}'}\varphi(x) dx+\left | \int
_{\mathbb{R}^{N}}(u^{2}-v^{2})\varphi dx \right |\\ &=2\int
_{\mathbb{R}^{N}}(u-v)^{m(p-1)}\mathcal{M}_{{s}'}\varphi(x) dx+\left | \int
_{\mathbb{R}^{N}}(u^{2}-v^{2})\varphi dx \right |, \end{align*} where
${s}'=sp/2$ and in $(a)$ we have used that $(u^{m}-v^{m})\leqslant (u-v)^{m}$.
By \eqref{380}, we have $$\left | \int _{\mathbb{R}^{N}}(u^{2}-v^{2})\varphi dx
\right |\leq L(2)\int _{\mathbb{R}^{N}} \left |u-v\right | dx=L(2)\int
_{\mathbb{R}^{N}}(u-v)\varphi dx.$$ 
Since the exponent $m(p-1)$ lies between
$0$ and $1$, we get by Hölder's inequality
\begin{align*}
    &\int _{\mathbb{R}^{N}}(u-v)^{m(p-1)}\mathcal{M}_{sp/2}\varphi(x) dx\\
&\leq \left ( \int _{\mathbb{R}^{N}}(u(x,t)-v(x,t))\varphi (x)dx \right
)^{m(p-1)}\\
&\quad\left ( \int _{\mathbb{R}^{N}}\frac{\left |
\mathcal{M}_{sp/2}\varphi(x) \right |^{\frac{1}{1-m(p-1)}}}{\varphi
(x)^{\frac{m(p-1)}{1-m(p-1)}}}dx \right )^{1-m(p-1)}. \end{align*} 
In the last
expression, we choose $\varphi \in \mathcal{C}(s,p,m)$ and we have called
$C(\varphi)$ is finite. We can write \begin{align*}
   & \left | _{0}^{C}\textrm{D}_{t}^{\alpha }\int _{\mathbb{R}^{N}}(u(x,t)-v(x,t))\varphi (x)dx \right|\\
   &\leqslant\left ( \int _{\mathbb{R}^{N}}(u(x,t)-v(x,t))\varphi (x)dx \right )^{m(p-1)} C(\varphi)^{1-m(p-1)}\\
   &\quad+L(2)\int _{\mathbb{R}^{N}}(u-v)\varphi dx.
\end{align*} By \ref{mkmk}, we can get $$ \left |
_{0}^{C}\textrm{D}_{t}^{\alpha }X(t) \right|\leqslant C(\varphi
)^{1-m(p-1)}X(t)^{m(p-1)}+L(2)X(t).$$ From Young's inequality, let $a=C(\varphi
)^{1-m(p-1)},b=X(t)^{m(p-1)},q=\frac{1}{m(p-1)}>1,p=\frac{1}{1-m(p-1)}$, then
we have \begin{align*}
    &C(\varphi )^{1-m(p-1)}X(t)^{m(p-1)}\\
    &\leqslant \varepsilon ^{\frac{1}{1-m(p-1)}}C(\varphi )+\frac{m(p-1)}{\varepsilon ^{\frac{1}{m(p-1)}}}X(t).
\end{align*} So we obtain \begin{equation}\label{sdsde}
     _{0}^{C}\textrm{D}_{t}^{\alpha }X(t)\leqslant \varepsilon ^{\frac{1}{1-m(p-1)}}C(\varphi )+[\frac{m(p-1)}{\varepsilon ^{\frac{1}{m(p-1)}}}+L(2)]X(t).
\end{equation} and by Lemma \ref{lemma2.14} and $a^{\alpha }-b^{\alpha
}\leqslant (a-b)^{\alpha },0<\alpha<1$ ,the above fractional differential
inequality \ref{sdsde} on $(t_{1},t_{2})$ with $t_{1},t_{2}\geq 0$ gives the
result \ref{dfdf}. 
\end{pro}
\begin{remark}
    The main references for the proof of this lemma are \cite{vazquez2022growing} and \cite{bonforte2014quantitative}.The first difference from above two literature is that we replace the partial derivative $\partial _{t}u$ in the equation with the Caputo fractional derivative $ _{0}^{C}\textrm{D}_{t}^{\alpha }u$ and the nonlinear fractional  diffusion term $(-\Delta )_{p}^{s}u^{m}$. Secondly, there are no partial terms in the above two literature, and the main reference \cite{borikhanov2022qualitative} for the treatment of non-local terms is scaled down. Finally, we also take advantage of the scaling of fractional differential inequalities in the proof process.
\end{remark}
\subsection{Proof of Theorem \ref{oplkhtrfv}}\label{28uijybb}
\begin{definition}\label{love}
 A function $u$ is a weak solution to the problem \eqref{100} if:
\begin{itemize}
    \item $u\in L^{2}((0,T];H_{0}^{1}(\mathbb{R}^{N}))\cap C((0,T];L^{1}(Q_{T}))$, where $Q_{T}=\mathbb{R}^{N}\times (0,T)$ and $$\left | u \right |^{m-1}u\in L_{loc}^{2}((0,T];W^{s,p}(\mathbb{R}^{N}));$$
    \item identity
    \begin{align*}
        &\int_{0}^{T}\int _{\mathbb{R}^{N}}u\frac{\partial^{\alpha} \varphi }{\partial t^{\alpha}}dxdt\\
        & +\int_{0}^{T}\iint_{\mathbb{R}^{2N}} \frac{(u^{m}(x,t)-u^{m}(y,t))^{p-1}(\varphi(x,t)-\varphi (y,t))}{\left | x-y \right |^{N+sp}}dxdydt\\
        &=\int_{0}^{T}\int _{\mathbb{R}^{N}}f\varphi dxdt
    \end{align*}
holds for every $\varphi \in C^{\infty}(\mathbb{R}^{N}\times [0,T])$;
    \item $u(\cdot,0)=u_{0}\geq0$ and $\int u_{0}(x)\varphi (x)dx<\infty $ for some admissible test function $\varphi \geq 0$ in the class $C(s,p,m)$.
\end{itemize} \end{definition}

\begin{definition}\label{sddd}
 We say that the  problem \eqref{100} with initial data $u_{0}\in L^{2}(\mathbb{R}^{N})$ has a unique strong solution $u\in C([0,T);L^{2}(\mathbb{R}^{N}))$, if moreover $$\partial _{t }^{\alpha }\in L^{\infty }((\tau,\infty);L^{1}(\mathbb{R}^{N}))\cap L^{2}(\mathbb{R}^{N}), (-\Delta)_{p}^{s}u \in L^{2}(\mathbb{R}^{N})$$ for every $t,\tau>0$.
\end{definition}
 \begin{remark}
 The above two definitions are mainly combined with references \cite{de2012general,tatar2017inverse,vazquez2021fractional}. On the one hand, \cite{tatar2017inverse,de2012general} gives the definition of strong solution for time-space fractional nonlinear diffusion equations with $\partial _{t }^{\alpha }\in L^{\infty }((\tau,\infty);L^{1}(\mathbb{R}^{N})),\tau>0$. On the other hand, \cite{vazquez2021fractional} gives the definition of  existence of strong solution for fractional $p$-Laplacian evolution equation with $u_{t}$ and $(-\Delta)_{p}^{s}u \in L^{2}(\mathbb{R}^{N})$ for every $t>0$.
 \end{remark}

 \begin{thm}\label{oplkhtrfv}
[Existence] Let $sp<1,0<m<\frac{1}{p-1}$ and let us consider the problem
\ref{100} exists a weak solution which is defined in $Q_{\infty}$ from
Definition \ref{love} with locally integrable initial data $u_{0}\geq0$. This
solution is continuous in the weighted space, $u\in
C([0,T];L^{1}(\mathbb{R}^{N},\varphi dx))$.
 \end{thm}
\begin{pro}
 \textbf{(1).} We follow the outline of proof of Theorem 3.1 of \cite{bonforte2014quantitative} and \cite{vazquez2022growing}, but some important changes are needed. Firstly,  by \cite{tatar2017inverse}, we denote the space $L^{2}(0,T;H_{0}^{1}(\mathbb{R}^{N}))$ by $V$ and define $Lu:=\left \langle \hat{L}u,\varphi \right \rangle,\hat{L}u:=\frac{\partial^{\alpha }u}{\partial t^{\alpha }},\hat{L}:D(\hat{L})\subset V\rightarrow V^{*}$ with the domain $$D(\hat{L})=\left \{ u\in V: \frac{\partial^{\alpha}u }{\partial t^{\alpha}}\in V^{*} \right \}.$$ We note that the operator $L$ is linear, densely defined and m-accretive, see \cite{bazhlekova2012existence}. Next, by \cite{vazquez2022growing}, let $\varphi \in \mathcal{C}$ and $\varphi _{R}$ the scaling of $\varphi$. Let $0\leq u_{0,n}\in L^{1}(\mathbb{R}^{N})\cap L^{\infty}(\mathbb{R}^{N})$ be a non-decreasing sequence of initial data $u_{0,n-1}\leq u_{0,n}$, converging monotonically to $u_{0}\in L^{1}(\mathbb{R}^{N},\varphi dx)$. By the Monotone Convergence Theorem, it follows that $\int _{\mathbb{R}^{N}}(u_{0}-u_{n,0})\varphi dx\rightarrow 0$ as $n\rightarrow \infty $.

 \textbf{(2).} We prove existence of the monotone limit of  the approximating solutions. From Definition \ref{sddd}, considering the unique strong solution $u_{n}(x,t)$ of Eq.\ref{100} with initial data $u_{0,n}$. By the comparison results of \cite{vazquez2021fractional}, we know that the sequence of solutions is a monotone sequence. The weighted estimates of previous section imply that the sequence is bounded in $L^{1}(\mathbb{R}^{N},\varphi dx)$ uniformly in $t\in [0,T].$
 \begin{align}
  \nonumber  & \left ( \int _{\mathbb{R}^{N}}u_{n}(x,t)\varphi (x)dx \right )^{1-m(p-1)}\\
 \nonumber   &\leqslant \left (  \int _{\mathbb{R}^{N}}u_{n}(x,0)\varphi (x)dx\right )^{1-m(p-1)}+K(\varphi )t^{\alpha [1-m(p-1)]}\\
&\leqslant \left (  \int _{\mathbb{R}^{N}}u_{0}(x)\varphi (x)dx\right
)^{1-m(p-1)}+K(\varphi )t^{\alpha [1-m(p-1)]}.
 \end{align}
  By the Monotone Convergence Theorem in $L^{1}(\mathbb{R}^{N},\varphi dx)$, we know that the solution $u_{n}(x,t)$ converge monotonically as $n\rightarrow \infty$ to a function $u(x,t)\in L^{\infty}((0,T);L^{1}(\mathbb{R}^{N},\varphi dx))$. We also have
  \begin{align*}
      \left ( \int _{\mathbb{R}^{N}}u(x,t)\varphi (x)dx \right )^{1-m(p-1)}&\leqslant\left (  \int _{\mathbb{R}^{N}}u_{0}(x)\varphi (x)dx\right )^{1-m(p-1)}\\
      &\qquad+K(\varphi )t^{\alpha [1-m(p-1)]}.
  \end{align*}

   \textbf{(3).}Next, we show that the function $u(x,t)$ is a weak solution to Eq. \ref{100} in $\mathbb{R}^{N}\times [0,T]$. We know that each $u_{n}$ is a bounded strong solution according to Definition \ref{sddd} since the initial data $u_{0}\in L^{1}(\mathbb{R}^{N})\cap L^{\infty }(\mathbb{R}^{N})$. Therefore, for all $\psi \in C_{c}^{\infty }(\mathbb{R}^{N}\times [0,T])$ we have
   \begin{align*}
        &\int_{0}^{T}\int _{\mathbb{R}^{N}}u_{n}(x,t)\frac{\partial^{\alpha} \psi }{\partial t^{\alpha}}dxdt\\
        &=\int_{0}^{T}\iint_{\mathbb{R}^{2N}} \frac{(u_{n}^{m}(x,t)-u_{n}^{m}(y,t))^{p-1}(\psi (x,t)-\psi  (y,t))}{\left | x-y \right |^{N+sp}} dxdydt\\
        &\qquad+\int_{0}^{T}\int _{\mathbb{R}^{N}}f(x,t)\psi(x,t) dxdt
    \end{align*}
Taking the limit $n\rightarrow \infty $ in the first line is easy:
\begin{equation*}
    \lim_{n\rightarrow \infty }\int_{0}^{T}\int _{\mathbb{R}^{N}}u_{n}(x,t)\frac{\partial^{\alpha} \psi }{\partial t^{\alpha}}dxdt=\int_{0}^{T}\int _{\mathbb{R}^{N}}u(x,t)\frac{\partial^{\alpha} \psi }{\partial t^{\alpha}}dxdt
\end{equation*} since $\psi$ is compactly supported and we already know that
$u_{n}(x,t)\rightarrow u(x,t)$ in $L_{loc}^{1}$.

\textbf{(4).}On the other hand, the integral in the second line $I(u_{n})$, is
well defined and can be estimated uniformly. We argue as before, by using the
inequality \begin{equation*}
    \left | u_{n}^{m}(x,t)-u_{n}^{m}(y,t)  \right |^{p-1}\leq u_{n}(x,t)^{m(p-1)}+u_{n}(y,t)^{m(p-1)},
\end{equation*} and bounding the first of the two ensuring integrals by
$$I_{1}(u_{n})=\int_{\mathbb{R}^{N}} dt\left ( \int
_{\mathbb{R}^{N}}u_{n}(x,t)^{m(p-1)}M(x,t)dx \right ).$$ where
$$M(x,t)=\int_{\mathbb{R}^{N}}\frac{\left | \psi(x,t)-\psi (y,t) \right
|}{\left | x-y \right |^{N+sp}}dy.$$ Due to the regularity of $\psi$, the last
integral is bounded above by some $\overline{M}_{\psi }(x)$ that behaves like
$C(1+\left | x \right |)^{-(N+sp)}$ independently of $t$. Hence,
\begin{align*}
    I_{1}(u_{n})&\leq \int_{\mathbb{R}^{N}} dt\left ( \int
_{\mathbb{R}^{N}}u_{n}(x,t)\varphi dx \right )^{m(p-1)}\\
&\quad\left ( \int
_{\mathbb{R}^{N}}\frac{\overline{M}_{\psi }(x)^{\frac{1}{1-m(p-1)}}}{\varphi
^{\frac{m(p-1)}{1-m(p-1)}}}dx \right )^{1-m(p-1)}
\end{align*}
so that finally
$$I_{1}(u_{n})\leq C(\psi,\varphi )\int _{\mathbb{R}^{N}}dt\left ( \int
_{\mathbb{R}^{N}}u_{n}(x,t)\varphi dx \right )^{m(p-1)}\leq CT$$ The second
integral $I_{2}(u_{n})$ is treated similarly by exchanging $x$ and $y$.

\textbf{(5).} We call that $u\geq u_{n}$, and we may apply the argument of the
previous paragraph to $v_{n}=u-u_{n}$ as follows $$I(u)-I(u_{n})=\int
_{0}^{T}\iint_{\mathbb{R}^{2N}}\frac{A(x,y,t)(\psi(x,t)-\psi(y,t) )}{\left |
x-y \right |^{N+sp}}dxdydt$$ with
$$A(x,y,t)=(u^{m}(x,t)-u^{m}(y,t))^{p-1}-(u_{n}^{m}(x,t)-u_{n}^{m}(y,t))^{p-1}.$$
Using the numerical inequality $\left | a^{m}-b^{m} \right |\leq 2^{1-m}\left |
a-b \right |^{m}$ and $(a)$ we conclude that for all possible values of
$(u(x,t)-u(y,t))$ and $(u_{n}(x,t)-u_{n}(y,t))$ we have \begin{align*}
    &\left | A(x,y,t) \right |\\
    &\leq \left | (u(x,t)-u(y,t))^{m(p-1)}-(u_{n}(x,t)-u_{n}(y,t))^{m(p-1)} \right |\\
&\leq c(m,p)\left | (u(x,t)-u(y,t))-(u_{n}(x,t)-u_{n}(y,t)) \right |^{m(p-1)}\\
&\leq c(m,p)\left | v_{n}(x,t)-v_{n}(y,t) \right |^{m(p-1)}, \end{align*} and
then we continue as before to prove that
\begin{align*}
    &\left | I_{1}(u)-I_{1}(u_{n}) \right |\\
&\leq C\left ( \iint(u-u_{n})\varphi dxdt \right )^{m(p-1)}\\
&\qquad\left (
\iint\frac{M(x,t)^{\frac{1}{1-m(p-1)}}}{\varphi^{\frac{m(p-1)}{1-m(p-1)}} }dxdt
\right )^{1-m(p-1)}, \end{align*} so that by virtue of the previous estimates
on the weighted convergence of $u_{n}(x,t)\rightarrow u(x,t)$, then we obtain
$$\left | I(u)-I(u_{n}) \right |\rightarrow 0\quad \text{as}\quad n\rightarrow
\infty .$$

\textbf{(6).}The solutions constructed above for $0\leq u_{0}\in
L^{1}(\mathbb{R}^{N},\varphi dx)$ satisfy the weighted estimates \ref{dfdf}, so
that
\begin{align}
    &\left | \int _{\mathbb{R}^{N}}u(x,t)\varphi (x)dx-\int _{\mathbb{R}^{N}}u(x,\tau )\varphi (x)dx \right |\\
    &\leq \frac{C_{1}(2\varepsilon)^{\frac{1}{1-m(p-1)}}}{\alpha \Gamma (\alpha )}\left | t-\tau  \right |^{\alpha }
\end{align} 
which gives the continuity in $L^{1}(\mathbb{R}^{N},\varphi
dx)$. Therefore, the initial trace of this solution is given by $u_{0}\in
L^{1}(\mathbb{R}^{N},\varphi dx)$. \end{pro}
 \begin{remark}
    The proof of the above solution mainly refers to the proof process outline of the Theorem 3.2 of \cite{bonforte2014quantitative}\cite{vazquez2022growing}. Through the above two papers, the uniqueness of the solution for the problem \eqref{100} mainly proves that solution $v\leq u$, and then the same argument can prove $u\leq v$.
\end{remark}

\end{appendices}

\bibliography{aipsamp}% Produces the bibliography via BibTeX.

\end{document}